\documentclass[12pt]{amsart}
\usepackage{oldlfont}
\newcommand{\HH}{{\cal H}}
\newcommand{\K}{{\cal K}}
\newcommand{\E}{{\cal E}}

\newcommand{\J}{{\bf J}}
\newcommand{\U}{\Upsilon}
\newcommand{\sm}{\setminus}
\newcommand{\al}{\alpha}
\newcommand{\eps}{\varepsilon}

\newcommand{\Z}{{\Bbb Z}} \newcommand{\R}{{\Bbb R}} \newcommand{\C}{{\Bbb C}}
\author{V.~A.~Vassiliev}
\address {Independent Moscow University, B.~Vlas'evskii per., 11.
Moscow, 121002, RUSSIA}
\email{vassil@vassil.mccme.rssi.ru}
\thanks{Supported in part by RFBR
(project 95-01-00846a)
and Netherlands Organization for Scientific Research
(NWO), project 47.03.005. Research at MSRI supported in part by
NSF grant DMS-9022140.}

\title{Topology of two-connected graphs and
homology of spaces of knots}

\begin{document}

\begin{abstract}
We propose a new method of computing cohomology groups
of spaces of knots in $\R^n$, $n \ge 3$,
based on the topology of configuration spaces and two-connected graphs,
and calculate all such classes of order $\le 3.$
As a byproduct we define the higher indices, which invariants of
knots in $\R^3$ define at arbitrary singular knots. More generally,
for any finite-order cohomology class of the space of knots
we define its principal symbol, which lies in a cohomology group of a certain
finite-dimensional configuration space and characterizes our class
modulo the classes of smaller filtration.
\end{abstract}

\maketitle

\section{Introduction}

The knots, i.e. smooth embeddings $S^1 \to \R^n$, $n \ge 3,$
form an open dense subset in the space ${\cal K} \equiv C^\infty(S^1,\R^n)$.
Its complement $\Sigma$ is the {\it discriminant set},
consisting of maps, having selfintersections or singularities.
Any cohomology class $\gamma \in H^i(\K \setminus \Sigma)$
of the space of knots can be described as the linking number
with an appropriate chain of codimension $i+1$ in $\K$ lying in $\Sigma$.

In \cite{V1} a method of constructing some of such (co)homology classes
(in particular, knot invariants) was proposed. For $n=3$ the 0-dimensional
classes, arising from this construction, are exactly the ``finite-type
knot invariants'', and for $n>3$ this method provides
a complete calculation of all cohomology groups of
knot spaces in $\R^n$.

However, the precise calculations by this method are very complicated.
The strongest results obtained by now are as follows.

D.~Bar-Natan calculated the $\C$-valued knot invariants of orders
$\le 9$ (and used for this several weeks of computer work
and the Kontsevich's realization theorem), see \cite{BN}.

T.~Stanford wrote a program, realizing the algorithm
of \cite{V1} for computing $\Z$-in\-va\-ri\-ants, and found
all such invariants of orders $\le 7$.

D.~Teiblum and V.~Turchin (also using a computer) found the first non-trivial
1-dimensional cohomology class of the space of knots in $\R^3$
(which is of order 3)
and proved that there are no other positive-dimensional cohomology classes
of order 3.
For a description of this class, see \S \ 4.4.2.

By the natural periodicity of cohomology groups of spaces of knots in $\R^n$
with different $n$
(similar to the fact that all algebras $H^*(\Omega S^m)$ with $m$
of the same parity are isomorphic up to a scaling of dimensions)
these results can de extended to cohomology of spaces of knots in
$\R^n$ with arbitrary odd $n$: to any knot invariant of order $i$ in
$\R^3$ there corresponds a $(n-3)i$-dimensional cohomology class for any odd
$n$, and the Teiblum--Turchin class is the origin of a series of
$((n-3)i+1)$-dimensional cohomology classes in spaces of knots in $\R^n$, $n$
odd.

Below I describe another, in a sense opposite method
of calculating these cohomology groups, based on a different
filtration of the discriminant variety.
On the level of knot invariants in $\R^3$,
this method is more or less equivalent to the Bar-Natan's
calculus of Chinese Character Diagrams and gives a partial explanation
of its geometrical sense.
Using this method, it is possible to repeat ``by hands'' the calculation of
Teiblum--Turchin (which was nontrivial even for computer)
and to obtain a similar result for the case of even $n$.

The idea of reducing the cohomology of spaces of nonsingular geometrical
objects (such as e.g. polynomials of a fixed degree in $\C^1$ or $\R^1$
without multiple roots) to the (Alexander dual) homology of the
complementary discriminant set was proposed by V.~I.~Arnold in
\cite{A1}, see also \cite{A2}. This idea proved to be very fruitful,
because the discriminant is a naturally stratified set (whose open
strata consist of ``equisingular'' objects). Many homology classes
of the discriminant can be calculated with the help of its
filtration, defined by this stratification.

The ``reversed'' filtrations of (some natural resolutions of) the discriminant
sets, introduced in \cite{V0}, simplify very much these calculations
and allow to solve similar problems in the multidimensional case,
see e.g. \cite{VAMS}, \cite{V4}.
Moreover, the spectral sequences, defined by such filtrations, are functorial
with respect to the inclusion of (finite-dimensional) functional
spaces, and thus give rise to stable sequences, calculating cohomology
of complements of discriminants in infinite-dimensional spaces like the
space of knots or of $C^\infty$-functions without complicated singularities.

This reversion is a continuous generalization of the combinatorial
formula of inclusions and exclusions: if we wish to calculate
the cardinality of a finite union of finite sets, then, instead of
counting the points in ``open strata'' (i.e., for any $k,$ in all $k$-fold
intersections of these sets, from which all the $(k+1)$-fold intersections
are removed) we can to count separately all points in all sets, then to
distract from the
obtained sum the correction term, corresponding to all 2-fold intersections
(not taking into account the fact that some of these double intersections
points are also triple points), then to add the correction terms, corresponding
to all triple points, etc. All results of \cite{V0}--\cite{V5}, concerning the
cohomology of complements of discriminants, were obtained by this method.

However, in the calculation of the first term $E_1$ of the spectral
sequence from \cite{V1}, converging to the cohomology of knot spaces (in
particular to knot invariants), we used an auxiliary spectral sequence, based on
the ``natural'' sub-stratification of the resolved discriminant.  The main
portion of hard calculations in the theory of finite-order knot invariants is
the calculation in this auxiliary sequence (especially in its terms responsible
for the 0-dimensional cohomology).  In this paper we reverse also this
filtration.  From the point of view of the new method the previous calculations
are just the computation of cellular homology of certain configuration spaces,
which can be studied by more ``theoretical'' methods.

As a byproduct, we obtain the notion of the {\it index}, which any
knot invariant assigns to any (finitely-degenerate) singular knot in $\R^3$.
In the standard
theory of finite-order invariants
this index was considered for the simplest points of $\Sigma$,
i.e. for immersions having only transverse double points or
(in a more implicit way) at most one triple point.

In the general situation this function is not numerical: it takes
values in certain homology groups related with singular knots.
Say, for an immersion having one generic $k$-fold
self-intersection point this group is the $(2k-4)$-th homology
group of the complex of {\it two-connected
graphs\footnote{Unlike the topological terminology, in the combinatorics
a graph with $k$ nodes is called two-connected if it is connected
(i.e. joins any two of these $k$ nodes), and,
moreover, removing from it any node together
with all incident edges,
we obtain also a connected graph.}} with $k$ nodes,
in particular, according to \cite{BBLSW} and \cite{T},
is $(k-2)!$-dimensional. In fact, the information provided by the
index of such a singular knot
is equivalent to the ``totality of all extensions of our invariant
to all Chinese Character Diagrams of order $k-1$ with exactly $k$ legs'',
see \cite{BN2}.

More generally, for any finite-order cohomology class of the space of knots
we define its {\it principal symbol}, which lies in a cohomology group of a
certain finite-dimensional configuration space and characterizes our class
modulo the classes of smaller filtration.

Our first calculations lead to some essential problems
in the homological
combinatorics and representation theory, see \cite{BBLSW}, \cite{T},
and the first answers indicate the
existence of a rich algebraic structure behind it, which probably will
allow one to guess many invariants and cohomology classes
of arbitrary orders.
\medskip

{\sc Notation.} For any topological space $X,$
$\bar H_*(X)$ denotes the {\it Borel--Moore homology group} of $X,$
i.e. the homology group of its one-point compactification reduced modulo the
added point.
\medskip

{\bf Acknowledgments.} I thank very much A.~Vaintrob,
who about January 1994
asked me whether it is possible
to define generalized indices  of complicated singular knots, and
A.~Bj\"or\-ner, V.~Turchin and V.~Welker for their interest and work
on combinatorial aspects of the problem.

\section{Complexes of connected and two-connected
graphs}

We consider only graphs without loops and multiple edges, but
maybe with isolated nodes.

A graph with $k$ nodes is {\it connected,} if any two of its nodes are
joined by a chain of its edges. A graph is
$l$-{\it connected} if it is connected,
and removing from it any $l-1$ nodes together with all incident
edges, we obtain again a connected graph (with $k-l+1$ nodes).

The set of all graphs with given $k$ nodes generates an (acyclic)
simplicial complex. Namely, consider the simplex $\Delta(k)$,  whose
vertices are in the one-to-one correspondence with all
$\binom{k}{2}$ edges of the complete graph with these $k$ nodes.
The faces of this simplex are subgraphs
of the complete graph: indeed, any face is
characterized by the set of its vertices,
i.e. by a collection of edges of the complete graph.
All faces of the simplex $\Delta(k)$ form an acyclic simplicial complex, which
also will be denoted by $\Delta(k)$. A generator of this complex is a
graph with ordered edges, while permuting the edges we send such a
generator to $\pm $ itself depending on the parity of the permutation.
We will always choose the
generator, corresponding to the lexicographic order of
edges induced by some fixed order of initial $k$ nodes.
The boundary of a graph is a formal sum of all
graphs obtained from it by removing one of its edges, taken with
coefficient 1 or $-1$.

Faces, corresponding to non-connected graphs, form a subcomplex of
the complex $\Delta(k)$. The corresponding quotient complex
is denoted by $\Delta^1(k)$ and is called the {\it complex of connected graphs}.
In a similar way, the complex $\Delta^l(k)$ of $l$-{\it connected graphs} is
the quotient complex of $\Delta(k)$, generated by all faces, corresponding
to $l$-connected graphs.
\medskip

{\sc Theorem 1.} a) For any $k,$ the group $H_i(\Delta^1(k))$
is trivial for all
$i \ne k-2,$ and $H_{k-2}(\Delta^1(k)) \simeq \Z^{(k-1)!}.$

b) A basis in the group $H_{k-2}(\Delta^1(k))$ consists of all linear
(homeomorphic to a segment) graphs, one of whose endpoints is fixed.
\medskip

Statement a) of this theorem is a corollary of a theorem of
Folkman \cite{Fo}, see e.g. \cite{B}. A proof of b) and another
proof of a), based on the
Goresky--MacPherson formula for the homology of plane arrangements,
is given in \cite{V2}.
\medskip

{\sc Theorem 2.} For any $k,$ the group $H_i(\Delta^2(k))$ is trivial
if $i \ne 2k-4$ and is isomorphic to $\Z^{(k-2)!}$ if $i=2k-4.$ \medskip
\medskip

This theorem was proved independently and almost simultaneously by
Eric Babson, Anders Bj\"orner, Svante Linusson, John Shareshian,
and Volkmar Welker, on one hand, and by Victor Turchin on the other,
see \cite{BBLSW}, \cite{T}.
\medskip

{\sc Example 1.} Suppose that $k=3$ and the original nodes are numbered by
1, 2 and 3. The simplex $\Delta(3)$ is a triangle, whose vertices are called
$(1,2),$ $(1,3)$ and $(2,3).$  Among its 7 faces only four correspond to
connected graphs, namely, all faces of dimension 1 or 2. In particular, the
homology group $H_i(\Delta^1(3))$  is trivial if $i \ne 1$ and is isomorphic to
$\Z^2$ if $i=1.$

The unique face of $\Delta(3)$, corresponding to a 2-connected graph,
is the triangle itself. Thus $H_i(\Delta^2(3))$ is trivial if $i \ne 2$
and is isomorphic to $\Z$ if $i=2.$
\medskip

{\sc Example 2.}
The simplex $\Delta(4)$ has $\binom{4}{2} =6$ vertices. The corresponding
complex of 2-connected graphs consists of the simplex itself,
all 6 its faces of dimension 4, and 3 faces of dimension 3, corresponding
to all cycles of length 4. It is easy to calculate that
$H_i(\Delta^2(4))=0$ for $i \ne 4$ and $H_4(\Delta^2(4)) \simeq \Z^2.$ Namely,
this homology group is generated by three $4$-chains, any of which
is the difference of two graphs
with 5 edges, obtained from the complete graph by removing edges,
connecting complementary pairs of points,
see Fig.~\ref{bchains}.
\begin{figure}
\begin{center}
\unitlength=1.00mm
\special{em:linewidth 0.4pt}
\linethickness{0.4pt}
\begin{picture}(28.00,7.00)(2,2)
\put(3.00,6.00){\makebox(0,0)[cc]{1}}
\put(3.00,2.00){\makebox(0,0)[cc]{2}}
\put(13.00,2.00){\makebox(0,0)[cc]{3}}
\put(13.00,6.00){\makebox(0,0)[cc]{4}}
\put(29.00,6.00){\makebox(0,0)[cc]{4}}
\put(29.00,2.00){\makebox(0,0)[cc]{3}}
\put(19.00,2.00){\makebox(0,0)[cc]{2}}
\put(19.00,6.00){\makebox(0,0)[cc]{1}}
\put(16.00,4.00){\makebox(0,0)[cc]{{\large $-$}}}
\put(5.00,1.00){\line(1,0){6.00}}
\put(11.00,1.00){\line(0,1){6.00}}
\put(11.00,7.00){\line(-1,0){6.00}}
\put(5.00,7.00){\line(0,-1){6.00}}
\put(5.00,1.00){\line(1,1){6.00}}
\put(21.00,7.00){\line(0,-1){6.00}}
\put(21.00,1.00){\line(1,0){6.00}}
\put(27.00,1.00){\line(0,1){6.00}}
\put(27.00,7.00){\line(-1,0){6.00}}
\put(21.00,7.00){\line(1,-1){6.00}}
\end{picture}
\qquad
\unitlength=1.00mm
\special{em:linewidth 0.4pt}
\linethickness{0.4pt}
\begin{picture}(28.00,7.00)(2,2)
\put(3.00,6.00){\makebox(0,0)[cc]{1}}
\put(3.00,2.00){\makebox(0,0)[cc]{2}}
\put(13.00,2.00){\makebox(0,0)[cc]{3}}
\put(13.00,6.00){\makebox(0,0)[cc]{4}}
\put(29.00,6.00){\makebox(0,0)[cc]{4}}
\put(29.00,2.00){\makebox(0,0)[cc]{3}}
\put(19.00,2.00){\makebox(0,0)[cc]{2}}
\put(19.00,6.00){\makebox(0,0)[cc]{1}}
\put(16.00,4.00){\makebox(0,0)[cc]{{\large $-$}}}
\put(5.00,1.00){\line(0,1){6.00}}
\put(5.00,7.00){\line(1,-1){6.00}}
\put(11.00,7.00){\line(-1,0){6.00}}
\put(5.00,1.00){\line(1,1){6.00}}
\put(11.00,7.00){\line(0,-1){6.00}}
\put(21.00,7.00){\line(0,-1){6.00}}
\put(21.00,1.00){\line(1,1){6.00}}
\put(27.00,1.00){\line(-1,0){6.00}}
\put(21.00,7.00){\line(1,-1){6.00}}
\put(27.00,1.00){\line(0,1){6.00}}
\end{picture}
\qquad
\unitlength=1.00mm
\special{em:linewidth 0.4pt}
\linethickness{0.4pt}
\begin{picture}(28.00,7.00)(2,2)
\put(3.00,6.00){\makebox(0,0)[cc]{1}}
\put(3.00,2.00){\makebox(0,0)[cc]{2}}
\put(13.00,2.00){\makebox(0,0)[cc]{3}}
\put(13.00,6.00){\makebox(0,0)[cc]{4}}
\put(29.00,6.00){\makebox(0,0)[cc]{4}}
\put(29.00,2.00){\makebox(0,0)[cc]{3}}
\put(19.00,2.00){\makebox(0,0)[cc]{2}}
\put(19.00,6.00){\makebox(0,0)[cc]{1}}
\put(16.00,4.00){\makebox(0,0)[cc]{{\large $-$}}}
\put(11.00,1.00){\line(-1,0){6.00}}
\put(5.00,1.00){\line(1,1){6.00}}
\put(5.00,7.00){\line(0,-1){6.00}}
\put(11.00,1.00){\line(-1,1){6.00}}
\put(5.00,7.00){\line(1,0){6.00}}
\put(27.00,7.00){\line(0,-1){6.00}}
\put(21.00,1.00){\line(1,0){6.00}}
\put(27.00,1.00){\line(-1,1){6.00}}
\put(21.00,7.00){\line(1,0){6.00}}
\put(27.00,7.00){\line(-1,-1){6.00}}
\end{picture}
\end{center}
\caption{Basic chains for 2-connected graphs with 4 vertices}
\label{bchains}
\end{figure}

Such basic chains are numbered by unordered partitions of four points into
two pairs and satisfy one relation: the sum of all three chains is equal
to the boundary of the complete graph.

\section{Simplest examples of indices of singular knots}

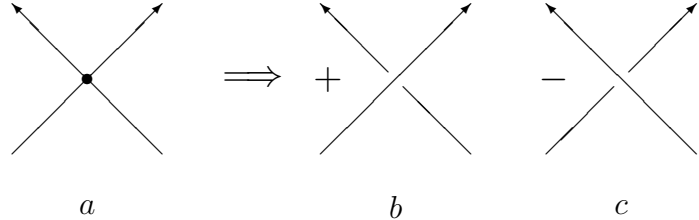
\begin{figure}
\unitlength=1.00mm
\special{em:linewidth 0.4pt}
\linethickness{0.4pt}
\begin{picture}(101.00,30.00)
\put(10.00,10.00){\vector(1,1){20.00}}
\put(30.00,10.00){\vector(-1,1){20.00}}
\put(20.00,20.00){\circle*{1.33}}
\put(51.00,10.00){\vector(1,1){20.00}}
\put(71.00,10.00){\line(-1,1){9.00}}
\put(60.00,21.00){\vector(-1,1){9.00}}
\put(101.00,10.00){\vector(-1,1){20.00}}
\put(92.00,21.00){\vector(1,1){9.00}}
\put(90.00,19.00){\line(-1,-1){9.00}}
\put(91.00,3.00){\makebox(0,0)[cc]{$c$}}
\put(61.00,3.00){\makebox(0,0)[cc]{$b$}}
\put(20.00,3.00){\makebox(0,0)[cc]{$a$}}
\put(52.00,20.00){\makebox(0,0)[cc]{{\large $+$}}}
\put(82.00,20.00){\makebox(0,0)[cc]{{\large $-$}}}
\put(42.00,20.00){\makebox(0,0)[cc]{{\large $\Longrightarrow$}}}
\end{picture}
\caption{Local perturbations of a self-intersection point}
\label{si}
\end{figure}

For any immersion $\phi: S^1 \to \R^3$ with exactly $i$ transverse crossings
we can consider all $2^i$ possible small
resolutions of this singular immersion,
replacing any of its self-intersection points
(see Fig.~\ref{si}a) by either undercrossing
or overcrossing. Any such local resolution can be invariantly
called positive (see Fig.~\ref{si}b)) or negative
(Fig.~\ref{si}c)), see \cite{V1}.
The sign of the
entire resolution is defined as $(-1)^{\# -},$ where $\# -$ is the number
of {\it negative} local resolutions in it.

Given a knot invariant, the index (or ``$i$-th jump'')
of our singular immersion is
defined as the sum of values of this invariant
at all positive resolutions minus
the similar sum over negative resolutions. Invariants of order
$\le j$ can be defined as those taking zero index at all immersions
with $ >j$ transverse crossings.

\begin{figure}
\unitlength=1.00mm
\special{em:linewidth 0.4pt}
\linethickness{0.4pt}
\begin{picture}(124.00,85.00)
\put(50.00,37.00){\vector(1,0){5.00}}
\put(50.00,37.00){\vector(0,1){5.00}}
\put(50.00,37.00){\vector(-1,-1){4.00}}
\put(55.00,34.00){\makebox(0,0)[cc]{$x$}}
\put(53.00,42.00){\makebox(0,0)[cc]{$z$}}
\put(49.00,31.00){\makebox(0,0)[cc]{$y$}}
\put(28.00,45.00){\vector(1,0){15.00}}
\put(38.00,40.00){\vector(0,1){16.00}}
\put(43.00,50.00){\vector(-1,-1){10.00}}
\put(38.00,45.00){\circle*{1.33}}
\put(23.00,12.00){\circle*{1.33}}
\put(28.00,17.00){\circle*{1.33}}
\put(55.00,17.00){\circle*{1.33}}
\put(55.00,14.00){\circle*{1.33}}
\put(75.00,39.00){\circle*{1.33}}
\put(71.00,39.00){\circle*{1.33}}
\put(62.00,71.00){\circle*{1.33}}
\put(57.00,66.00){\circle*{1.33}}
\put(30.00,75.00){\circle*{1.33}}
\put(30.00,70.00){\circle*{1.33}}
\put(9.00,47.00){\circle*{1.33}}
\put(13.00,47.00){\circle*{1.33}}
\put(77.00,18.00){\makebox(0,0)[cc]{$3$}}
\put(39.00,4.00){\makebox(0,0)[cc]{$4$}}
\put(5.00,25.00){\makebox(0,0)[cc]{$5$}}
\put(6.00,68.00){\makebox(0,0)[cc]{$6$}}
\put(45.00,82.00){\makebox(0,0)[cc]{$1$}}
\put(79.00,61.00){\makebox(0,0)[cc]{$2$}}
\put(93.00,27.00){\circle*{1.33}}
\put(105.00,27.00){\circle*{1.33}}
\put(119.00,27.00){\circle*{1.33}}
\put(115.00,27.00){\circle*{1.33}}
\put(104.00,27.00){\oval(22.00,4.00)[t]}
\put(112.00,27.00){\oval(14.00,4.00)[b]}
\put(118.00,37.00){\circle*{1.33}}
\put(107.00,37.00){\circle*{1.33}}
\put(93.00,37.00){\circle*{1.33}}
\put(97.00,37.00){\circle*{1.33}}
\put(105.50,37.00){\oval(25.00,4.00)[t]}
\put(102.00,37.00){\oval(10.00,4.00)[b]}
\put(93.00,47.00){\circle*{1.33}}
\put(104.00,47.00){\circle*{1.33}}
\put(108.00,47.00){\circle*{1.33}}
\put(119.00,47.00){\circle*{1.33}}
\put(98.50,47.00){\oval(11.00,4.00)[t]}
\put(113.50,47.00){\oval(11.00,4.00)[b]}
\put(119.00,57.00){\circle*{1.33}}
\put(115.00,57.00){\circle*{1.33}}
\put(105.00,57.00){\circle*{1.33}}
\put(93.00,57.00){\circle*{1.33}}
\put(106.00,57.00){\oval(26.00,4.00)[t]}
\put(110.00,57.00){\oval(10.00,4.00)[b]}
\put(119.00,67.00){\circle*{1.33}}
\put(107.00,67.00){\circle*{1.33}}
\put(93.00,67.00){\circle*{1.33}}
\put(97.00,67.00){\circle*{1.33}}
\put(100.00,67.00){\oval(14.00,4.00)[t]}
\put(108.00,67.00){\oval(22.00,4.00)[b]}
\put(119.00,77.00){\circle*{1.33}}
\put(108.00,77.00){\circle*{1.33}}
\put(104.00,77.00){\circle*{1.33}}
\put(93.00,77.00){\circle*{1.33}}
\put(100.50,77.00){\oval(15.00,4.00)[t]}
\put(111.50,77.00){\oval(15.00,4.00)[b]}
\put(124.00,77.00){\makebox(0,0)[cc]{$1$}}
\put(124.00,67.00){\makebox(0,0)[cc]{$2$}}
\put(124.00,57.00){\makebox(0,0)[cc]{$3$}}
\put(124.00,47.00){\makebox(0,0)[cc]{$4$}}
\put(124.00,37.00){\makebox(0,0)[cc]{$5$}}
\put(124.00,27.00){\makebox(0,0)[cc]{$6$}}
\put(23.00,7.00){\vector(0,1){16.00}}
\put(55.00,8.00){\vector(0,1){15.00}}
\put(75.00,34.00){\vector(0,1){15.00}}
\put(62.00,67.00){\vector(0,1){10.00}}
\put(30.00,64.00){\vector(0,1){16.00}}
\put(9.00,44.00){\vector(0,1){9.00}}
\put(42.00,43.00){\oval(84.00,84.00)[]}
\put(42.00,43.00){\oval(34.00,34.00)[]}
\put(59.00,53.00){\line(5,3){25.00}}
\put(59.00,33.00){\line(5,-3){25.00}}
\put(25.00,33.00){\line(-5,-3){25.00}}
\put(25.00,53.00){\line(-5,3){25.00}}
\put(2.00,47.00){\line(1,0){16.00}}
\put(16.00,50.00){\line(-1,-1){11.00}}
\put(9.00,37.00){\line(0,1){5.00}}
\put(16.00,17.00){\line(1,0){6.00}}
\put(24.00,17.00){\line(1,0){9.00}}
\put(31.00,20.00){\line(-1,-1){12.00}}
\put(50.00,9.00){\line(1,1){7.00}}
\put(59.00,18.00){\line(1,1){5.00}}
\put(64.00,17.00){\line(-1,0){15.00}}
\put(65.00,39.00){\line(1,0){15.00}}
\put(67.00,35.00){\line(1,1){7.00}}
\put(76.00,44.00){\line(1,1){5.00}}
\put(66.00,75.00){\line(-1,-1){13.00}}
\put(52.00,66.00){\line(1,0){15.00}}
\put(62.00,77.00){\line(0,-1){10.00}}
\put(62.00,65.00){\line(0,-1){5.00}}
\put(34.00,79.00){\line(-1,-1){13.00}}
\put(24.00,70.00){\line(-1,0){5.00}}
\put(26.00,70.00){\line(1,0){9.00}}
\put(121.00,27.00){\line(-1,0){30.00}}
\put(91.00,37.00){\line(1,0){30.00}}
\put(121.00,47.00){\line(-1,0){30.00}}
\put(91.00,57.00){\line(1,0){30.00}}
\put(121.00,77.00){\line(-1,0){30.00}}
\put(91.00,67.00){\line(1,0){30.00}}
\put(42.00,1.00){\line(0,1){25.00}}
\put(42.00,60.00){\line(0,1){25.00}}
\put(35.50,40.00){\oval(5.00,10.00)[b]}
\put(43.00,47.50){\oval(6.00,5.00)[r]}
\end{picture}
\caption{Splittings of a triple point}
\label{split}
\end{figure}
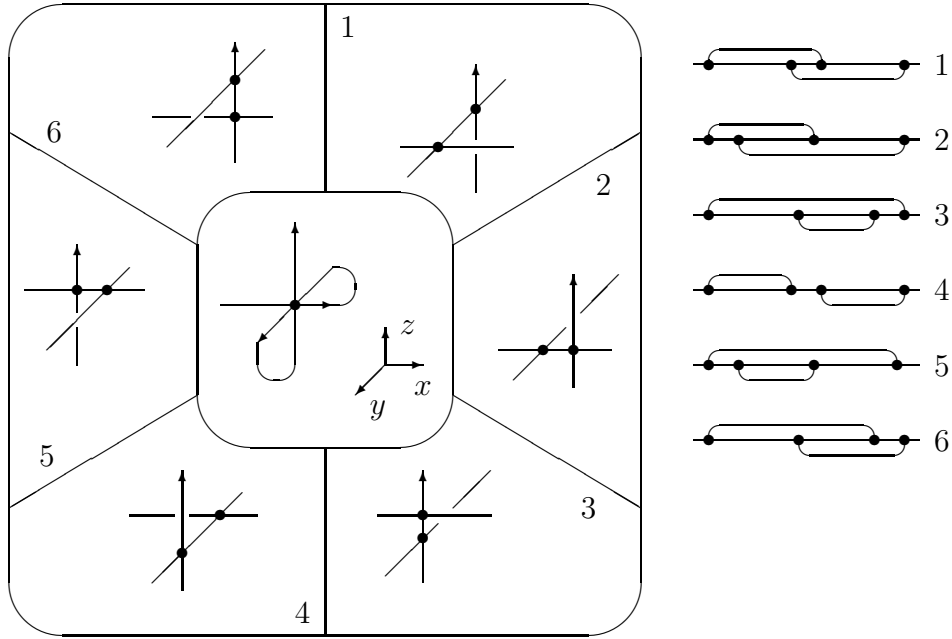

Further, let $\phi: S^1 \to \R^3$ be an immersion with $l$ transverse
double self-in\-ter\-sec\-ti\-ons and one generic triple point. There are 6
different perturbations of this triple point, splitting it into a
pair of double points, see Fig.~\ref{split} (= Fig.~15 in \cite{V1}). Thus we
get 6 immersions with $l+2$ double crossings each, and, given a knot invariant,
6 their indices $I(1), \ldots, I(6)$ of order $l+2.$ These indices are
dependent: they satisfy the {\it four-term relations}
\begin{equation}
\label{fterm}
I(1)-I(4)=I(2)-I(5)=I(6)-I(3).
\end{equation}

The common value of these three differences is a characteristic
assigned by our invariant to the initial singular immersion
with a triple point and
also is called its index, see \cite{V1}.

An equivalent definition of the order of invariants is as follows:
an invariant is of order $j$ if all
such indices assigned by it to all immersions with $\ge j-1$ double
and one triple generic self-intersection points vanish.

In \S \ 5.3 we define similar indices for maps
$S^1 \to \R^3$ having an arbitrary finite number of multiple
self-intersections or singular points.

\section{The discriminant set and its resolution}

In this section we recall basic facts from \cite{V1} concerning the
topological structure of the discriminant set in the space
$\K$ of all smooth maps $S^1 \to \R^n,$ $n \ge 3$.

For simplicity we mention this space as a
space of very large but finite dimension $\omega$.
A partial justification of this assumption uses
finite-dimensional approximations of $\K$, see \cite{V1}.
Below we indicate by the single quotes `, ' the non-rigorous assertions
using this assumption and needing a reference to \cite{V1} for such an
justification.

\subsection{Simplicial resolution of the discriminant.}

The {\it resolution} $\sigma$ of the discriminant set $\Sigma$ is constructed
as follows. Denote by $\Psi$ the space of all unordered pairs $(x,y)$ of
points of $S^1$ (allowing $x=y$); it is easy to see that $\Psi$ is
diffeomorphic to the closed M\"obius band. Consider a generic embedding
$I$ of $\Psi$ in a space $\R^N$ of ({\it very})$^2$
large dimension $N \gg \omega^2.$
For any discriminant map $\phi: S^1 \to \R^n$ consider all points
$(x,y) \subset \Psi$ such that either $x \ne y$ and $\phi(x) = \phi(y),$
or $x=y$ and $\phi'(x)=0.$ Denote by $\Delta(\phi)$ the convex hull in $\R^N$
of images of all such points under the embedding $I$. If $I$ is generic, then
this convex hull is a simplex, whose vertices coincide with all these images.
The space $\sigma$ is defined as the union of all simplices of the form
$\phi \times \Delta(\phi) \subset \K \times \R^N.$  The restriction on
$\sigma$ of the obvious projection $\K \times \R^N \to \K$ is proper and
induces a `homotopy equivalence' $\pi: \bar \sigma \to \bar \Sigma$
of one-point compactifications  of spaces $\sigma$ and $\Sigma.$
By the `Alexander duality',
the homology groups $\bar H_*(\sigma) \equiv \bar H_*(\Sigma)$
of these compactifications
`coincide' (up to a change of dimensions) with the cohomology
groups of the space of knots:
\begin{equation}
H^i(\K \setminus \Sigma) \simeq \bar H_{\omega-i-1} (\Sigma) \equiv
\bar H_{\omega-i-1}(\sigma).
\label{alex}
\end{equation}

The calculation of these groups is based on the natural stratification of
$\Sigma$ and $\sigma$.
Let $A$ be a non-ordered finite collection of natural
numbers, $A=(a_1, a_2, \ldots, a_{\#A}),$ any of which is not less than
$2,$ and
$b$ a non-negative integer. Set $|A|=a_1+ \cdots +a_{\#A}.$
An $(A,b)$-{\it configuration} is a collection of $|A|$ distinct points
in $S^1$ separated into groups of cardinalities $a_1, \ldots, a_{\#A},$
plus a collection of $b$ distinct points in $S^1$ (some of which can
coincide with the above $|A|$ points). For brevity,
$(A,0)$-configurations are called simply $A$-configurations.
A map $\phi: S^1 \to \R^n$ {\it respects}
an $(A,b)$-configuration if it glues together all points inside
any of its groups of cardinalities $a_1,\ldots, a_{\# A},$ and its derivative
$\phi'$ is equal to 0 at all the $b$ points of this
configuration.  For any $(A,b)$-configuration
the set of all maps, respecting it,
is an affine subspace in $\K$ of codimension $n(|A|-\#A+b)$;
the number $|A|-\#A+b$ is called the {\it complexity} of the configuration.
Two $(A,b)$-configurations are {\it equivalent} if they
can be transformed into one another by an orientation-preserving homeomorphism
$S^1 \to S^1.$

For any $(A,b)$-configuration $J$, $A=(a_1, \ldots, a_{\# A}),$ the
corresponding simplex $\Delta(J) \subset \R^N$ is defined as the simplex
$\Delta(\phi)$ for an arbitrary {\it generic}
map $\phi$ respecting $J;$ it has exactly
$\binom{a_1}{2} + \cdots +\binom{a_{\# A}}{2} + b$ vertices.

For any class ${\bf J}$ of equivalent $(A,b)$-configurations, the
corresponding ${\bf J}$-{\it block} in $\sigma$ is defined as the union of all
pairs $\{\phi,x\}$, where $\phi$ is a map $S^1 \to \R^n$, respecting some
$(A,b)$-configuration $J$ of this equivalence class, and $x$ is
a point of the simplex $\Delta(J)$.

\subsection{The main filtration in the resolved discriminant.}

The {\it main filtration} in $\sigma$ is defined as follows: its
term $\sigma_i$ is the union of all ${\bf J}$-blocks over all equivalence
classes $\J$ of $(A,b)$-configurations of complexities $\le i.$

This filtration is the unique useful filtration in entire
$\sigma$ and will not be revised. Indeed,
its term $\sigma_i \setminus \sigma_{i-1}$ is the space of an affine
bundle of dimension $\omega - ni$ over a finite-dimensional base,
which can easily be described: the one-point compactification of this base is a
finite cell complex, see \cite{V1}. The spectral sequence, `calculating' the
groups (\ref{alex}) and induced by this filtration, satisfies the condition

\begin{equation}
\label{ineq}
E_{p,q}^1=0  \mbox{ for } p(n-2)+q > \omega-1.
\end{equation}
(By the definition of the spectral
sequence, $E_{p,q}^1 = \bar H_{p+q}(\sigma_i \sm \sigma_{i-1}).)$

The Alexander dual {\it cohomological} spectral sequence
$E_r^{p,q}$ is obtained from this one by the formal change of indices,
\begin{equation}
\label{dual}
E_r^{p,q} \equiv E^r_{-p, \omega-1-q};
\end{equation}
it `converges' to some subgroups of groups $H^{p+q}(\K \sm \Sigma)$,
and the support of its term $E_1$ belongs to the wedge
\begin{equation}
\label{ineq2}
p \le 0, p(n-2)+q \ge 0.
\end{equation}

\subsubsection{Kontsevich's realization theorem.}

M.~Kontsevich \cite{K2} proved that this spectral sequence (over $\C$)
degenerates at the first term: $E_1^{p,q} \equiv E^{p,q}_\infty.$
At least if $n=3,$ then for the groups $E^{p,q}$ with $p+q=0$,
which provide knot invariants, this follows from his integral realization
of finite-order invariants, described in \cite{K}, \cite{BN2}.
Moreover, for arbitrary $n$ and $p<0$ exactly the same construction
proves the degeneration of the ``leading term'' of the column
$E^{p,*}$:
\begin{equation}
\label{ineq3}
E_1^{p,p(2-n)}(\C) \simeq E_\infty^{p,p(2-n)}(\C).
\end{equation}
In several other
cases (in particular if $-p$ is sufficiently small with respect to $n$) this
follows also from dimensional reasons. However, generally for
greater values of $q$ the proof is more complicated.

\subsubsection{The main filtration and connected graphs.}
This filtration induces an infinite filtration in (some subgroup of)
the cohomology
group $H^*(\K \sm \Sigma) \simeq \bar H_{\omega-*-1}(\sigma):$
an element of this group has order $i$ if it can be defined
as the linking number with a cycle
lying in the term $\sigma_i.$ In particular, for $*=0$ we get
a filtration in the space of knot invariants; it is easy to see
that this filtration
coincides with the elementary characterization
of finite-order invariants given in \S \ 3.

Let $J$ be an $(A,b)$-configuration of complexity $i$.
By construction, the corresponding $\J$-block is a fiber
bundle, whose base is the space of $(A,b)$-configurations
equivalent to $J$, and the fiber over such a configuration
$J'$ is the direct product of an affine space of dimension
$\omega - ni$ (consisting of all maps respecting this configuration)
and the simplex $\Delta(J')$. Consider this $\J$-block
as the space of the fiber bundle, whose base is the
$(\omega-ni)$-dimensional affine bundle over the previous base, and the
fibers are simplices $\Delta(J')$.
Some points of these simplices belong not only to the
$i$-th term $\sigma_i$ of our filtration, but even to the $(i-1)$-th term.
These points form a simplicial subcomplex of the simplex $\Delta(J')$,
let us describe it. Any face
of $\Delta(J')$ can be depicted by $\# A$ graphs with $a_1, \ldots, a_{\# A}$
vertices and $b$ signs $\pm$. Indeed, to any vertex of $\Delta(J')$
there corresponds
either an edge, connecting two points inside some of our $\# A$ groups of
points, or one of last $b$ points of the $(A,b)$-configuration $J'$.
\medskip

{\sc Proposition 1} (see \cite{V1}).
{\it A face of
$\Delta(J')$ belongs to $\sigma_{i-1}$ if and only if either at least one of
corresponding $\# A$ graphs is non-connected or at least one of $b$ points does
not participate in its picture (= participates with sign $-$).} \quad $\Box$
\medskip

Denote by $\Delta^1(J')$ the quotient complex of
$\Delta(J')$ by the
union of all faces belonging to $\sigma_{i-1}$.
Theorem 1 implies immediately the following statement.
\medskip

{\sc Proposition 2.} {\it
The quotient complex $\Delta^1(J')$
is acyclic in all dimensions other
than $|A|-\# A+b-1 \equiv i-1$, while}
\begin{equation}
\bar H_{|A|-\# A+b-1}(\Delta^1(J')) \simeq
\otimes_{j=1}^{\# A} \Z^{(a_j-1)!}.
\end{equation}

\subsection{Non-compact knots}

Simultaneously with the usual knots,
i.e., embeddings
$S^1 \to \R^n,$ we will consider {\it non-compact knots}, i.e. embeddings
$\R^1 \to \R^n$ coinciding with the standard linear embedding outside
some compact subset in $\R^1$. The space of all smooth maps with this behavior
at infinity will be denoted by $K$, and the
discriminant $\Sigma \subset K$ again is
defined as the set of all such maps,
having singularities and self-intersections.
In this case the configuration space $\Psi$, participating in the
construction of the resolution $\sigma$,
is not the M\"obius band, but the closed half-plane
$\R^2 /\{(t,t')\equiv(t',t)\},$ which we usually will realize as the
half-plane $\{(t,t')|t \le t'\}.$
\medskip

There is an obvious one-to one correspondence between isotopy classes
of standard and non-compact knots in $\R^3$, and spaces of invariants,
provided by the above spectral sequences
in both theories, naturally coincide. On the other hand, the
CW-structure on resolved discriminants in spaces of non-compact
knots is easier; this is the reason why in \cite{V1}
only the non-compact knots were considered. In the next
section 5 we develop some new techniques for calculating
cohomology groups of spaces of compact or non-compact knots in $\R^n$,
and in section 6 (respectively, 7) apply it to the calculation of cohomology
classes of order $\le 3$ of the space of non-compact knots in
$\R^n$ (respectively, to the classes of order $\le 2$ of the
space of compact knots).

\subsection{Ancient auxiliary filtration in the term
$\sigma_i \sm \sigma_{i-1}$ of the main filtration.}

By definition, the term $E^1_{i,q}$ of the main spectral sequence is
isomorphic to $\bar H_{i+q}(\sigma_i \sm \sigma_{i-1}).$
To calculate these groups,
the {\it auxiliary filtration} in the space $\sigma_i \sm \sigma_{i-1}$
was defined in \cite{V1}. Namely, its term $G_\alpha$ is the union
of all $\J$-blocks, such that
$\J$ is an equivalence class of $(A,b)$-configurations of complexity
$i$, consisting of $\le \alpha$ geometrically distinct points.
In particular, $G_{2i}= \sigma_i \sm \sigma_{i-1}$ and
$G_{i-1}=\emptyset$.
By Proposition 2, $\bar H_j(G_\al)=0$ if $j \ge \omega-i(n-1)+\al$.
\medskip

\subsubsection{Example: knot invariants in $\R^3$.} If $n=3$, then
$\bar H_j(G_{2i-2})=0$
for $j \ge \omega-2,$ in particular, calculating the group
$\bar H_{\omega-1}(\sigma_i \sm \sigma_{i-1})$,
we can ignore all $\J$-blocks with auxiliary filtration $\le 2i-2:$
\begin{equation}
\bar H_{\omega-1}(\sigma_i \sm \sigma_{i-1}) \simeq
\bar H_{\omega-1}((\sigma_i \sm \sigma_{i-1}) \sm G_{2i-2}).
\label{red}
\end{equation}

The remaining part
$(\sigma_i \sm \sigma_{i-1}) \sm G_{2i-2}$ consists of
$\J$-blocks satisfying the following condition.
\medskip

{\sc Definition.} An $(A,b)$-configuration $J$
(and the corresponding $\J$-block) of complexity $i$
is {\it simple}, if one of three
is satisfied:\footnote{We denote by $2^{\times m}$ the expression
$2, \ldots, 2$ with $2$ repeated $m$ times}

I. $A=(2^{\times i})$, $b=0$;

II. $A=(2^{\times (i-1)})$, $b=1$ and the last point does not coincide
with any of $2i-2$ points, defining the $A$-part of the configuration;

III. $A=(3,2^{\times (i-2)})$, $b=0.$
\medskip

The complex $\Delta^1(J)$, corresponding to a configuration $J$ of
type I or II,
consists of unique cell of dimension $i-1,$ i.e. of the simplex
$\Delta(J)$ itself, while for $A=(3, 2, \ldots, 2)$ it consists of one
$i$-dimensional simplex and certain three of its faces of dimension $i-1$.

A chain complex, calculating groups (\ref{red})
(and, moreover, all the groups
$\bar H_*(\sigma_i \sm \sigma_{i-1})$)
was written out in \cite{V1}; its part calculating the top-dimensional
group (\ref{red}) can be described in the following terms.

An  $(A,b)$-configuration of type I (or the equivalence
class of such configurations) can be depicted by a
{\it chord diagram}, i.e. a collection of $i$
chords spanning some $i$ pairs of distinct points of $S^1$ or $\R^1$.
$(A,b)$-configurations of type III are depicted in \cite{V1}
by similar diagrams
with $i+1$ chords, three of which form a triangle and remaining $i-2$
have no common endpoints.
In \cite{BN2} this configuration is depicted by the same collection
of $i-2$ chords and an {\bf Y}-wise star connecting the points of the triple
with a point not on the line (or circle).
Configurations of types I and III are called
$[i]$- and $\langle i \rangle$-configurations, respectively, see
\cite{V1}.

An element of the group (\ref{red}) is a linear combination of several
$\J$-blocks of type I and interior parts (swept out by open $i$-dimensional
simplices)
of blocks of type III. This linear combination should satisfy the
homology condition
close to all strata of codimension 1 in  $\sigma_i:$ any such stratum should
participate in the algebraic boundary of this combination with coefficient 0.
Close to such a stratum, corresponding to an $(A,b)$-configuration $J$
of type II, there is unique stratum of maximal dimension: this is the
stratum of type I, whose $(A,b)$-configuration can be obtained from that of
$J$ by replacing its unique singular point by a small chord,
connecting two points close to it. The corresponding homological
condition is as follows: any stratum of maximal dimension, whose
chord diagram can be
obtained in this way (i.e. has a chord, whose endpoints are not
separated by endpoints of other chords) should participate in the linear
combination with zero coefficient.

Further, any configuration of type III defines three strata of codimension 1
in $\sigma_i$ (some of which can coincide) swept out by three
$(i-1)$-dimensional faces of
corresponding simplices $\Delta(J').$ Such a stratum can be described by an
$\langle i \rangle$-configuration,
in which one chord, forming the triangle, is erased.
This stratum is incident to three strata of maximal dimension:
one stratum swept out by open $i$-dimensional simplices
from the same $\J$-block
of type III, defined by the same
$\langle i \rangle$-configuration,
and two strata of type I corresponding to chord diagrams, obtained from our
configuration by diversing the endpoints of our non-complete triangle:

\unitlength=1.00mm
\special{em:linewidth 0.4pt}
\linethickness{0.4pt}
\begin{picture}(130.00,10.00)
\put(15.50,5.00){\oval(13.00,4.00)[t]}
\put(29.00,5.00){\oval(14.00,4.00)[t]}
\put(60.50,5.00){\oval(13.00,4.00)[t]}
\put(78.00,5.00){\oval(14.00,4.00)[t]}
\put(91.00,5.00){\makebox(0,0)[cc]{$\&$}}
\put(43.00,5.00){\vector(1,0){6.00}}
\put(5.00,5.00){\line(1,0){35.00}}
\put(52.00,5.00){\line(1,0){35.00}}
\put(95.00,5.00){\line(1,0){35.00}}
\put(9.00,5.00){\circle*{1.33}}
\put(22.00,5.00){\circle*{1.33}}
\put(36.00,5.00){\circle*{1.33}}
\put(85.00,5.00){\circle*{1.33}}
\put(71.00,5.00){\circle*{1.33}}
\put(67.00,5.00){\circle*{1.33}}
\put(54.00,5.00){\circle*{1.33}}
\put(118.50,5.00){\oval(17.00,6.00)[t]}
\put(106.00,5.00){\oval(16.00,4.00)[t]}
\put(127.00,5.00){\circle*{1.33}}
\put(114.00,5.00){\circle*{1.33}}
\put(110.00,5.00){\circle*{1.33}}
\put(98.00,5.00){\circle*{1.33}}
\end{picture}

The homological condition is as follows: the difference of coefficients, with
which two latter strata participate in the linear combination, should be
equal to the coefficient of the stratum swept out by triangles.
In particular, all three such differences, corresponding to erasing
any of three edges of the chord triangle, should be equal to one another.
These {\it three} equalities (among which only two are independent)
are called the {\it 4-term relations,} cf. \S \ 3.

The common value of these three differences is a characteristic
of our $(3,2,\ldots,2)$-configuration (and of the element $\alpha$
of the homology
group (\ref{red})): it is equal to the coefficient,
with which the main stratum of
the $\J$-block of type III participates in the cycle $\alpha$. It is natural
to call it the {\it index} of this cycle at this configuration.

In a similar way, any element $\al$ of the group (\ref{red}) assigns an
index to any $i$-chord diagram: this is the coefficient, with which the
corresponding $\J$-block of type I
participates in the cycle realizing $\al$.  \medskip

The group (\ref{red}) is thus canonically isomorphic to the group
of all $\Z$-valued functions on the set of $i$-chord diagrams, which a)
take zero value on all diagrams having chords
not crossed by other chords of the diagram, and b) satisfy
4-term relations defined by all possible $\langle i \rangle$-configurations.

When $i$ grows, this system of equations grows exponentially,
and, which is even worse, the answers do not satisfy any
transparent rule, see \cite{BN}.

\subsubsection{The Teiblum--Turchin cocycle}

In this subsection we consider the resolution of the discriminant in the space
of non-compact knots.

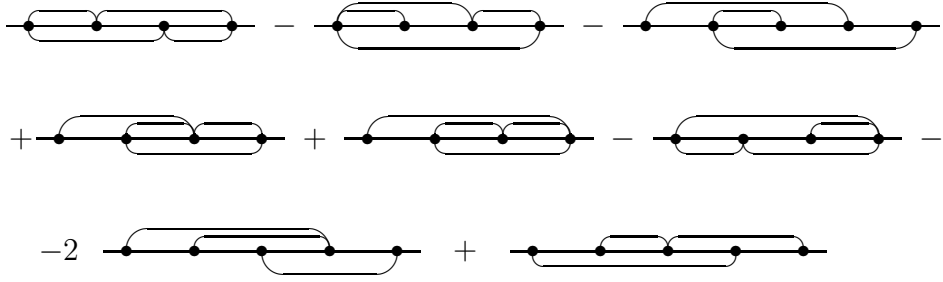
\begin{figure}
\unitlength 1mm
\special{em:linewidth 0.4pt}
\linethickness{0.4pt}
\begin{picture}(124,38)
\put(4,20){\line(1,0){33}}
\put(34,20){\circle*{1.33}}
\put(7,20){\circle*{1.33}}
\put(16,20){\circle*{1.33}}
\put(25,20){\circle*{1.33}}
\put(45,20){\line(1,0){33}}
\put(75,20){\circle*{1.33}}
\put(48,20){\circle*{1.33}}
\put(66,20){\circle*{1.33}}
\put(86,20){\line(1,0){33}}
\put(116,20){\circle*{1.33}}
\put(89,20){\circle*{1.33}}
\put(98,20){\circle*{1.33}}
\put(107,20){\circle*{1.33}}
\put(0,35){\line(1,0){33}}
\put(30,35){\circle*{1.33}}
\put(3,35){\circle*{1.33}}
\put(12,35){\circle*{1.33}}
\put(21,35){\circle*{1.33}}
\put(41,35){\line(1,0){33}}
\put(71,35){\circle*{1.33}}
\put(53,35){\circle*{1.33}}
\put(62,35){\circle*{1.33}}
\put(82,35){\line(1,0){42}}
\put(112,35){\circle*{1.33}}
\put(85,35){\circle*{1.33}}
\put(94,35){\circle*{1.33}}
\put(103,35){\circle*{1.33}}
\put(13,5){\line(1,0){42}}
\put(52,5){\circle*{1.33}}
\put(25,5){\circle*{1.33}}
\put(43,5){\circle*{1.33}}
\put(67,5){\line(1,0){42}}
\put(97,5){\circle*{1.33}}
\put(70,5){\circle*{1.33}}
\put(79,5){\circle*{1.33}}
\put(88,5){\circle*{1.33}}
\put(111.50,20){\oval(9,4)[t]}
\put(107,20){\oval(18,4)[b]}
\put(102.50,20){\oval(27,6)[t]}
\put(82,20){\makebox(0,0)[cc]{$-$}}
\put(70.50,20){\oval(9,4)[t]}
\put(66,20){\oval(18,4)[b]}
\put(41,20){\makebox(0,0)[cc]{$+$}}
\put(20.50,20){\oval(9,4)[t]}
\put(29.50,20){\oval(9,4)[t]}
\put(25,20){\oval(18,4)[b]}
\put(16,20){\oval(18,6)[t]}
\put(7,5){\makebox(0,0)[cc]{$-2$}}
\put(16,5){\circle*{1.33}}
\put(106,5){\circle*{1.33}}
\put(29.50,5){\oval(27,6)[t]}
\put(43,5){\oval(18,6)[b]}
\put(34,5){\circle*{1.33}}
\put(34,5){\oval(18,4)[t]}
\put(61,5){\makebox(0,0)[cc]{$+$}}
\put(83.50,5){\oval(9,4)[t]}
\put(97,5){\oval(18,4)[t]}
\put(83.50,5){\oval(27,4)[b]}
\put(2,20){\makebox(0,0)[cc]{$+$}}
\put(123,20){\makebox(0,0)[cc]{$-$}}
\put(7.50,35){\oval(9,4)[t]}
\put(12,35){\oval(18,4)[b]}
\put(25.50,35){\oval(9,4)[b]}
\put(21,35){\oval(18,4)[t]}
\put(37,35){\makebox(0,0)[cc]{$-$}}
\put(44,35){\circle*{1.33}}
\put(48.50,35){\oval(9,4)[t]}
\put(53,35){\oval(18,6)[t]}
\put(66.50,35){\oval(9,4)[t]}
\put(57.50,35){\oval(27,6)[b]}
\put(78,35){\makebox(0,0)[cc]{$-$}}
\put(98.50,35){\oval(27,6)[t]}
\put(98.50,35){\oval(9,4)[t]}
\put(107.50,35){\oval(27,6)[b]}
\put(121,35){\circle*{1.33}}
\put(57,20){\circle*{1.33}}
\put(61.50,20){\oval(9,4)[t]}
\put(61.50,20){\oval(27,6)[t]}
\put(93.50,20){\oval(9,4)[b]}
\end{picture}
\caption{The Teiblum--Turchin cocycle}
\label{tt}
\end{figure}

Let  $\J$ be an arbitrary equivalence class of
$(A,b)$-configurations of complexity $i$ in $\R^1$.
Denote by $\tilde B(\J)$ the intersection of the corresponding $\J$-block
with $\sigma_i \sm \sigma_{i-1}$. This set $\tilde B(\J)$
consists of several cells, which are in one-to-one
correspondence  with all possible collections of $\# A$ connected
graphs with $a_1, \ldots, a_{\# A}$ vertices respectively; for some
examples see Fig.~\ref{tt}. Such cells, corresponding to all equivalence
classes $\J$ of $(A,b)$-configurations of complexity $i$, define a cell
decomposition of the quotient space $\sigma_i / \sigma_{i-1}$.
\medskip

{\sc Theorem} (D.~Teiblum, V.~Turchin). {\it The linear combination of
suitably oriented (see \cite{V1}) cells of $\sigma_3 \setminus \sigma_2$
shown in Fig.~\ref{tt} defines a non-trivial $(\omega-2)$-dimensional
Borel--Moore homology class of $\sigma_3 \setminus \sigma_2$ (and hence also an
$(\omega-2-3(n-3))$-dimensional
homology class of the term $\sigma_3 \sm \sigma_2$ of the similar
resolution of the space of singular non-compact
knots in $\R^n$ for any odd $n \ge 3$).}
\medskip

On the other hand, the group $\bar H_{\omega-3-3(n-3)}(\sigma_2)$
is trivial (see
\cite{V1}), hence this class can be extended to a well--defined class in
$\bar H_{\omega-2-3(n-3)}(\sigma)$ and in $H^{1+3(n-3)}(K \setminus \Sigma)$.
If $n>3$, then by dimensional reasons this class is nontrivial, however
its nontriviality for $n=3$ and any explicit realization are unknown yet.

\section{New auxiliary filtration}

For any $(A,b)$-configuration $J$ denote by $\rho(J)$ the number of
{\it geometrically distinct} points in it. In particular, for any
$A$-configuration $J$, $\rho(J)=|A|$. The greatest possible value of
$\rho(J)$ over all $A$- or $(A,b)$-configurations of
complexity $i$ is equal to $2i.$
\medskip

{\sc Definition.} The {\it stickiness} of an $A$-configuration $J$ of complexity
$i$ is the number $2i-\rho(J)$.
The {\it reversed auxiliary filtration} $\Phi_0 \subset
\cdots \subset \Phi_{i-1}$ in the term $\sigma_i \sm \sigma_{i-1}$
of the main filtration is defined as follows:
its term $\Phi_\al$ is the closure of the union of
$\J$-blocks over all equivalence classes $\J$ of
$A$-configurations of stickiness $\le \al.$
\medskip

If $J$ is an $(A,b)$-configuration with $b>0$, then
the corresponding $\J$-block belongs to the closure of an
$\tilde \J$-block, where $\tilde J$ is an
$\tilde A$-configuration, $\tilde A=A \cup 2^{\times b}$.
Thus $\Phi_{i-1} = \sigma_i \sm \sigma_{i-1}$.
\smallskip

So, to any finite-order cohomology class $v$ of the space of knots
in $\R^n$ (in particular, to any invariant of knots in $\R^3$)
there correspond two numbers:
the first, $i(v),$ is its order, and the second, $r(v),$
is its reversed filtration, i.e. the minimal reversed filtration in
$\sigma_{i(v)} \sm \sigma_{i(v)-1}$ of realizing it
cycles in $\sigma_{i(v)}.$

\subsection{First examples.} $\ $

{\sc Notation.} For any topological space $X$, the {\it $k$-th
configuration space} $B(X,k)$ is the space of all subsets of
cardinality $k$ in  $X.$
$\pm \Z$ is the local system of groups on the configuration
space $B(X,k)$, locally isomorphic to $\Z$ and such that elements of
$\pi_1(B(X,k)),$ defining odd permutations of $k$ points, act on its fibers
as multiplication by $-1$. The Borel--Moore homology group
$\bar H_*(B(\Psi,k), \pm \Z)$ is the homology group of locally
finite chains with coefficients in this local system, cf. \cite{V4}.
\medskip

In the following two examples we consider the discriminant in the space of
non-compact knots.
\medskip

{\sc Example 3.} If $i=1$, then the entire space
$\sigma_1 \equiv \Phi_0$ is the space of an
$(\omega-n)$-dimensional affine fiber bundle over
the manifold $\Psi,$ diffeomorphic to the closed half-plane.
In particular, all its Borel--Moore homology groups are trivial
and there are no non-trivial cohomology classes of order 1
for an arbitrary $n$.
\medskip

{\sc Example 4.} Our filtration of the set $\sigma_2 \sm \sigma_1$
consists of two terms $\Phi_0 \subset \Phi_1.$
$\Phi_0$ is the space of a fiber bundle, whose base is the
configuration space $B(\Psi,2)$, and the fiber is the product of
an open interval and an affine space of codimension $2n$
in the space $K$. The generator of the group
$\pi_1(B(\Psi,2)) \simeq \Z$ changes the orientation of the first
factor of the fiber and multiplies the orientation of the second
by $(-1)^n$. In particular,
the Borel--Moore homology group of the term $\Phi_0$ is isomorphic
(up to the shift of dimensions) to the group
$\bar H_*(B(\Psi,2))$ if $n$ is odd and to
$\bar H_*(B(\Psi,2), \pm \Z)$ if $n$ is even.
It is easy to calculate (see also Lemma 1 below)
that both these groups are trivial in all dimensions.

The term  $\Phi_1 \sm \Phi_0$ is the space of a fiber bundle,
whose base is the configuration space $B(\R^1,3)$
of triples of points in the line, and the fiber is the product
of an affine subspace of codimension $2n$ in $K$ and
the interior part of a triangle. The Borel--Moore
homology of this term is obviously isomorphic to $\Z$
in the dimension $\omega -2n+5$ and is trivial in all other dimensions.
Thus the group of cohomology classes of order 2 is isomorphic to
$\Z$ in dimension $2n-6$ and is trivial in all other dimensions.
\smallskip

{\it A comparison.} Calculating the
homology of $\sigma_2 \sm \sigma_1$
with the help of the ancient auxiliary filtration, we assign the
edges of these triangles to the same term as their interior parts;
as a consequence, they separate the space
$B(\Psi,2)$ into several cells, corresponding to combinatorial
types of chord diagrams. In the new approach
all these cells are joined together and form
unique manifold $B(\Psi,2)$ with a simple topology.

In a similar way, for any $i$ the space
$\sigma_i \sm \sigma_{i-1}$ consists of $i$ terms
$\Phi_0 \subset \cdots \subset \Phi_{i-1}$
of the reversed filtration, and the ultimate term
$\Phi_{i-1} \sm \Phi_{i-2}$ is the space of a fiber bundle, whose
base is the configuration space $B(\R^1,i+1),$ and the fiber is the
direct product of an affine subspace of codimension
$ni$ in the space $K$ and an
$\binom{i+1}{2}$-vertex simplex, from which all faces belonging to
$\Phi_{i-2}$ are removed. These faces are exactly those corresponding
to not 2-connected graphs, see Theorem 3 below.

Similarly, to any immersion
$\R^1 \to \R^n$, having $j<\infty$ self-in\-ter\-sec\-ti\-on points of
finite multiplicities $a_1, \ldots, a_j$,
we associate the homology group of the corresponding component
of the term $\Phi_\al \sm \Phi_{\al-1}$ of our filtration in
$\sigma_i \sm \sigma_{i-1},$ where
$i=(\sum_{m=1}^j a_m-j),$ $\al= \sum_{m=1}^j a_m - 2j.$
This group is isomorphic (up to a shift of dimensions)
to the tensor product of homology groups of complexes of two-connected
graphs with $a_1, \ldots, a_j$ nodes.

\subsection{Structure of the reversed filtration.}

The set $\Phi_\al \sm \Phi_{\al-1}$ of the new filtration
in $\sigma_i \sm \sigma_{i-1}$ consists of several components,
numbered by unordered decompositions of the number $2i-\al$ into
$i-\al$ subsets of cardinalities $\ge 2$. To any such decomposition
\begin{equation}
\label{a}
A=(a_1, a_2, \ldots, a_{\# A}), \quad \# A =i-\al, \quad
\sum a_l=2i-\al, \quad a_l \ge 2,
\end{equation}
there corresponds the closure of all $\J$-blocks defined by all
equivalence classes $\J$ of $A$-configurations.
Let us describe this closure explicitly.
\medskip

{\sc Definition.}
An $A$-{\it collection} is a collection of $|A|$ points in $\R^1$
separated into $\# A$ groups of cardinalities $a_1, \ldots, a_{\# A}$,
such that
the points inside any group of cardinality $>2$ are pairwise distinct.

A map $\phi: \R^1 \to \R^n$ {\it respects} an $A$-collection, if
$\phi(x)=\phi(y)$ for any two different points $x, y $ of
the same group, and $\phi'(x)=0$ for any point $x$ such
that $(x,x)$ is a group of cardinality $a_i=2$ of this collection.

The set of all maps $\phi: \R^1 \to \R^n,$ respecting an $A$-collection,
is an affine subspace in $K$ of codimension
$\le n(|A| - \# A).$ Indeed, any group of cardinality $a_j$
gives us $n(a_j -1)$ independent restrictions on $\phi$
(although the union of all $n\sum (a_j-1) \equiv n(|A|-\# A)$ such
conditions can be dependent).
\medskip

An $A$-collection is an $A$-{\it set} if these conditions are independent,
i.e. the codimension of the space of maps $\phi: \R^1 \to \R^n$, respecting
this collection, is equal to $n(|A|-\# A).$

For example, $(2,2,2)$-sets are  any subsets of cardinality $3$ in
$\Psi$ not of the form $((x< y),(x<z),(y<z))$.
\medskip

More generally, given an $A$-collection,
consider the graph
with $|A|$ vertices, corresponding to its points, and
edges of following two types: all points inside any group are connected by
black edges, and geometrically coincident points of different groups
are connected by white edges. (If there is a degenerate
group of the form $(x,x)$ and the point $x$ belongs also to other
group(s), then it is depicted by two points connected by a black
edge, and only one of these points is connected by white edges with
corresponding points of other groups.)
\medskip

{\sc Proposition 3.} {\it Our $A$-collection is an $A$-set if and only if
any  cycle of this graph, having no repeating nodes,
contains no white edges.
In particular, two different groups cannot have two common points.}
\medskip

{\it Proof.} Suppose that we have such a cycle containing white edges.
Removing an arbitrary its vertex, incident to a white edge of the cycle,
from the corresponding group of the $A$-collection (say, from the group
of cardinality $a_l$) we obtain an $\tilde A$-collection,
$\tilde A=(a_1, \ldots, a_l-1, \ldots, a_{\# A})$, defining the same
subspace in $K.$ Thus the codimension of this subspace is
no greater than $n(|A|-\# A-1).$ Conversely, if we have no such cycles, then
there exists a node of our graph, not incident to white edges.
Removing the corresponding point from the $A$-collection
(and cancelling the group containing this point if it is a group
of cardinality 2) we obtain an $\tilde A$-collection with
$|\tilde A| - \# \tilde A = |A|-\# A-1.$ On the other hand, erasing this point
we loss exactly $n$ linear conditions on the corresponding subspace in $K,$
which are independent on the others. Thus the ``if'' part of Proposition 3
follows by induction. \quad $\Box$
\medskip

An $A$-set and $A'$-set are {\it related} if the sets of maps
$\R^1 \to \R^n$, respecting them, coincide. For example,
the $(3)$-set $((x<y<z))$ and $(2,2)$-set $((x,y), (x,z))$ are related.
There is a partial order on any class of related sets:
the set $\bar \U$ is a {\it completion} of the related set
$\U$ if $\bar \U$ is obtained from $\U$ by replacing some two intersecting
groups by their union or by a sequence of such operations.
\medskip

The space of all $A$-sets with given $A$ is denoted by $C(A).$
\medskip

An open dense subset in $C(A)$ consists of $A$-configurations,
see \S \ 4.1 above. E.g., such a subset in $C(2,2,2)$ consists
of all points $((x,x'), (y,y'), (z,z')) \in B(\Psi,3)$ such that none two of
six numbers $x, x', \ldots, z'$ coincide.

\subsubsection{The complex $\Lambda(\U)$.}

Given an $A$-set $\U$ in $\R^1,$ consider the simplex $\Delta(\U)$
with $\sum_{j=1}^{\# A} \binom{a_j}{2}$ vertices, corresponding
to all pairs $(x,y)\in \Psi$ such that $x$ and $y$ belong to
one group of $\U$. A face of this simplex can be encoded by a collection
of $\# A$ graphs, whose edges span the pairs corresponding to
vertices of this face.

Such a face is called 2-{\it connected}, if all these
graphs are 2-connected.
\medskip

We will identify the nodes of these graphs with corresponding points in
$\R^1.$ \medskip

Denote by $\Lambda(\U)$ the union of interior points of all
2-connected faces of our simplex.

Theorem 2 implies immediately the following statement.
\medskip

{\sc Corollary 1.} {\it 1. The group $\bar H_l(\Lambda(\U))$
is trivial for all $l \ne 2|A|-3\# A - 1 $  and is  isomorphic
to $\otimes_{j=1}^{\# A} \Z^{(a_j-2)!}$ for $l = 2|A|-3\# A-1;$
this isomorphism is defined canonically up to permutations of
factors $\Z^{(a_j-2)!}$, corresponding to coinciding numbers $a_j$.

2. For any decomposition of the set $\U$ into two non-intersecting
subsets $\U',$ $\U''$, any of which is the union of several groups
of $\U,$ we have
\begin{equation}
\label{join}
\bar H_{* - 1}(\Lambda(\U)) \simeq
\bar H_*(\Lambda(\U')) \otimes \bar H_*(\Lambda(\U'')),
\end{equation}
where the lower index $*-1$ denotes the shift of grading. \quad $\Box$}
\medskip

(In other words, if we introduce the notation
$$\HH(\U) \equiv \bar H_{2|A|-3\# A -1}(\Lambda(\U))$$
for any $A$ and any $A$-set $\U$, then
$\HH(\U) \simeq \HH(\U') \otimes \HH(\U'')$.)
\medskip

Namely, if we have two cycles
$\gamma' \subset \Lambda(\U') \subset \Delta(\U')$,
$\gamma'' \subset \Lambda(\U'') \subset \Delta(\U'')$,
then the corresponding {\it join} cycle
$\gamma' * \gamma'' \subset \Delta(\U) \equiv
\Delta(\U') * \Delta(U'')$ is realized as the union of all open
intervals, connecting the points of $\gamma'$ and $\gamma''$.
Similarly, if $\gamma^j,$ $j=1, \ldots, \#A,$ are cycles defining some
elements of groups $\bar H_*(\Lambda(\U_j)) \simeq \Z^{(a_j-2)!}$, then
their join $ \gamma^1 * \cdots * \gamma^{\#A}$ is swept out by open
$(\#A-1)$-dimensional simplices, whose vertices are the points
of these cycles $\gamma^j$. The corresponding homology map
$\otimes \bar H_*(\Lambda(\U_j)) \to \bar H_{*+\#A-1}(\Lambda(\U))$
is an isomorphism and depends on the choice of the orientation of
these simplices (= the ordering of their vertices or, equivalently,
of the groups $\U_j \subset \U$).
\medskip

Let $A$ be a multiindex (\ref{a}) of complexity $i,$
and $\U$ an $A$-set in $\R^1.$ The simplex $\Delta(\U)$ can be realized
as the simplex in $\R^N$, spanned by
all points $I(x,y)$, $(x,y) \in \Psi,$ such that $x$ and
$y$ belong to the same
group of $\U.$ For any map $\phi \in K$, respecting $\U,$
the simplex $\phi \times \Delta(\U) \subset K \times \R^N$ belongs
to the term $\sigma_i$ of the main filtration, and its intersection
with the space $\sigma_i \sm \sigma_{i-1}$ belongs to the term
$\Phi_\al$ of the reversed filtration of this space,
where $\al = 2i-\rho(\U).$
\medskip

{\sc Proposition 4.} {\it If $\bar \U$ is a completion of $\U$ and $\bar \U
\ne \U,$ then the simplex $\Delta(\U)$ is a not two-connected face of
$\Delta(\bar \U).$} \medskip

{\it Proof.} Consider the two-color graph of the $A$-set $\U$
mentioned in Proposition 3. The graph, representing the face
$\Delta(\U) \subset \Delta(\bar \U)$, is obtained from it by contracting
all white edges.
Since $\U \ne \bar \U,$ there is at least one white
edge connecting two points of the same group of $\bar \U$.
Removing from $\Delta(\U)$ the vertex, obtained from this edge,
we get a graph, which by Proposition 3 is not connected.
\quad $\Box$
\medskip

{\sc Theorem 3.} {\it For any multiindex $A$ of the form (\ref{a})
and any $A$-set $\U$, a
 point of the simplex $\phi \times \Delta(\U)$
belongs to $\Phi_\al \sm \Phi_{\al-1}$ if and only if it is an interior
point of a 2-connected face.}
\medskip

{\it Proof.}  First we prove the ``only if'' part.
Consider a face of the simplex $\Delta(\U)$ such that one of corresponding
graphs $g_1, \ldots,$ $ g_{\# A}$, say $g_m,$  is connected but not
two-connected. We need to prove that this face belongs to
$\Phi_{\al-1}$. Since $\Phi_{\al-1}$ is closed in $\Phi_{\al},$
it is sufficient to prove this property
for all $A$-sets $\U$ from an arbitrary
dense subset in $C(A).$ In particular we can assume that $\U$
is an $A$-configuration (see \S \ 4.1), i.e. all its $|A|$ points
are pairwise distinct.

Let $y \in \R^1$ be a node of our  graph $g_m$ such that removing
it we split $g_m$ into two non-empty graphs $g'_m, g''_m$ with
$a'_m$ and $a''_m$ nodes respectively, $a'_m+a''_m=a_m-1.$
Let $A'$ be the multiindex $(a_1, \ldots, a_{m-1}, a'_m+1,a''_m+1,a_{m+1},
\ldots, a_{\# A})$ and $\U'$ an $A'$-set, related to $\U$, whose groups of
cardinalities $a_1, \ldots, $ $a_{m-1}, $
$a_{m+1}, \ldots, a_{\# A}$ coincide
with these for $\U$, and the last two coincide with the sets of vertices of
$g'_m,$ $g''_m,$ augmented by $y.$ Then our face of the simplex $\phi \times
\Delta(\U) \simeq \Delta(\U)$ lies also in the simplex $\phi \times
\Delta(\U').$ This simplex belongs to the closure of the set swept out by
similar simplices $\phi' \times \Delta(\U'(\tau)),$ $\tau \in (0,\eps],$
where
$\U'(\tau)$ are $A'$-sets with
$\rho(\U'(\tau))> \rho(\U')=\rho(\U),$ namely,
some their $|A'|-2$ points coincide with these for
$\U',$ and only two points $y$ of
the $m$-th and $(m+1)$-th groups are replaced by $y-\eps$ and $y+\eps$
respectively. Thus our face belongs to a lower
term of the reversed filtration.

Conversely, for any multiindex $A$ of the form (\ref{a}),
the closure of the union of all
${\bf J}$-blocks over all $A$-configurations $J$ consists of simplices of the
form $\phi \times \Delta(\U)$, where $\U$ is an $A$-set, and $\phi$ a map
respecting it. By Proposition 4, if $\U$ is such an $A$-set and $\bar \U$
a completion of $\U,$ $\bar \U \ne \U,$ then such a simplex is a not
two-connected face of $\Delta(\bar \U).$ \quad $\Box$
\medskip

{\sc Corollary 2.} {\it The component of $\Phi_\al \sm \Phi_{\al-1},$
corresponding to the multiindex $A$ of the form (\ref{a}),
is the space of a fiber
bundle over the space of all $A$-sets, whose fiber
over the $A$-set $\U$ is the
direct product of the complex $\Lambda(\U)$ and the
affine subspace of codimension $\ ni \ $ in
$K,$ consisting of all maps $\phi$, respecting $\U.$} \quad $\Box$
\medskip

Denote this component by $S(A)$.
\medskip

{\sc Corollary 3.} {\it For any $A$, the group $\bar H_*(S(A))$
is trivial in all dimensions greater than $\omega -1-(n-3)i.$}
\medskip

This follows immediately from  Theorems 2 and 3 and implies
formulas (\ref{ineq}) and (\ref{ineq2}).

\subsection{Index of a knot invariant at a complicated singular knot.}

Let $A$ be a multiindex $(a_1, \dots, a_{\#A})$ of complexity $i$,
$J$ an $A$-configuration,
consisting of groups $J_1, \ldots, J_{\#A}$, and $\phi$ an immersion
$\R^1 \to \R^3$, respecting $J$ but not respecting more complicated
configurations. Define the group $\Xi_*(J)$ as the tensor product
$$ \otimes_{m=1}^{\#A} \bar H_*(\Lambda(J_m))$$
(where $\Lambda(J_m)$ is the complex of two-connected graphs, whose
nodes are identified with points of the group $J_m$).

By Theorem 2 this graded group
 has unique nontrivial term in dimension  $2|A|-4\#A$
and is isomorphic to $\otimes_{m=1}^{\#A} \Z^{(a_j-2)!}$ in this dimension.
In particular, there is an isomorphism
\begin{equation}
\label{order}
\Xi_{*}(J) \equiv \bar H_{*+\#A-1}(\Lambda(J)).
\end{equation}
By the {\it join} construction (see \S \ 5.2.1) this isomorphism is defined
almost canonically, only up to a sign, which can be specified by any ordering of
our groups $J_m$ and is different for orders of different parity.  \medskip

Any knot invariant $V$ of order $i$
defines an element of this group $\Xi_*(J)$,
the {\it index} of the singular knot $\phi$; this index
generalizes similar indices of
not very complicated singular knots considered in \S\S \ 3 and 4.4.1.
Indeed, consider any small open contractible neighborhood $U$ of our point
$\phi$ in the space of discriminant maps equisingular
to $\phi$ (i.e. respecting equivalent configurations).
The complete pre-image of $U$ in the space
$\sigma_i$ is homeomorphic to the direct product of $U$ and
the simplex $\Delta(\phi)$. Remove from this simplex
all faces, lying in lower terms of the main and reversed
auxiliary filtrations. Remaining domain will be homeomorphic to
$U \times \Lambda(J)$. The intersection of our invariant (i.e. of the
cycle $\gamma \in \bar H_{\omega-1}(\sigma_i) $ realizing it) with this
domain defines (via the K\"unneth isomorphism) a class
in $\bar H_{2|A|-3\#A-1}(\Lambda(J))$ and hence, by the isomorphism
(\ref{order}), also in the group $\Xi_*(J)$; this class is the desired
index $\langle V| \phi \rangle$.
\medskip

To make this definition correct, we need to specify the orientation of
$U$, participating in the construction of the K\"unneth isomorphism. This
orientation consists of a) an  orientation of the space $\J$ of
configurations equivalent to $J$,  and b) an orientation of the space
$\chi(J)$ of maps $\R^1 \to \R^3$ respecting $J$.
To define the first orientation, we order all groups $J_m$ of $J$ in an
arbitrary way and order points of any group by their increase
in $\R^1$. Thus we get an
order of all points of $J$ (first there go points of the first group in
increasing order, then of the second, etc.), and thus also the orientation of
the configuration space $\J$.

Further, we will suppose that an orientation of the functional space
$K$ is fixed, then orientations of the subspace $\chi(J)$ can be identified
with its coorientations, i.e. orientations of the normal bundle. Let $X, Y, Z$
be fixed coordinates in $\R^3$, so that the map $\phi$ consists of three real
functions $X(t), Y(t), Z(t)$. Then a coorientation of $\chi(J)$ can be specified
by the differential form $[d(X($the second point $t_{1,2}$ of the first group
$J_1 \subset J)- X($the first point $t_{1,1}$ of this group$)) \wedge
d(Y(t_{1,2})-Y(t_{1,1})) \wedge d(Z(t_{1,2})-Z(t_{1,1}))] \wedge \ldots
\wedge [d(X(t_{1,a_1})-X(t_{1,a_1-1})) \wedge d(Y(t_{1,a_1})-Y(t_{1,a_1-1}))
\wedge d(Z(t_{1,a_1})-Z(t_{1,a_1-1}))] \wedge ($the same for the second group
$J_2) \wedge \ldots \wedge ($ the same for $J_{\#A})$.  \medskip

{\sc Lemma 1.} {\it 1. The orientation of $U$, defined by this pair
of orientations of $\J$ and $\chi(J)$, depends on the choice of
the order of groups $J_m \subset J$, namely, an odd permutation of
these groups multiplies it by $-1$.

2. This orientation will be preserved if we change the increasing order
of points in an arbitrary group $J_m$ by means of any cyclic permutation.}
\medskip

The proof is elementary. \quad $\Box$
\medskip

The first statement of lemma implies that our index
$\langle V| \phi \rangle$ does not depend
on the choice of the order of groups $J_m \subset J$ (formally participating
in its construction). Indeed, changing this order by an odd permutation,
we multiply by $-1$ both the K\"unneth isomorphism and the isomorphism
(\ref{order}).

The second statement allows us to define a similar index
also in the case of compact knots $S^1 \to \R^3$.
\medskip

Moreover, suppose that our discriminant point  $ \phi$ is generic in its
stratum of complexity $i,$ i.e. it does not belong to the closure
of the set of maps, respecting configurations of complexities
$>i$. (Typical examples of non-generic points are immersions,
having self-tangency
points or triple points, at which three local branches are
complanar.) Then {\it any} knot invariant (of arbitrary order)
defines in the same way its generalized index, taking its value
in the same homology group. Indeed, in this case the complete intersection
of our neighborhood $U$ in entire
$\sigma\setminus \sigma_{i-1}$ (and not in
$\sigma_i\setminus \sigma_{i-1}$) has
the same structure of the direct product $U \times \Lambda(J)$.

For any class $\J$ of equivalent $A$-configurations, all groups
$\Xi_*(J)$, $J \in \J$, are canonically isomorphic to one another,
thus we can define the group $\Xi_*(\J)$ as any of them.
\medskip

{\sc Lemma 2.}
{\it Let $V$ be an  invariant of order $i$ of non-compact
knots in $\R^3$, and $\J$ a
class of equivalent $A$-{\it configurations} in $\R^1$ of complexity
$|A|-\# A=i.$
Then all indices $\langle V| \phi \rangle \in \Xi_*(\J)$ for all
generic immersions $\phi$, respecting $A$-configurations
$J \in \J$, coincide.} \quad $\Box$
\medskip

This (obvious) lemma allows us to define the index $\langle V| \J \rangle \in
\Xi_*(\J)$ as the common value of these indices $\langle V| \phi \rangle$.
\medskip

These indices satisfy the natural homological condition
(the analog of the $STU$-relation of \cite{BN2}). Let $J,J'$ be
two $A$-configurations from neighboring equivalence classes, i.e. $J'$ can
be obtained from $J$ by the permutation of exactly two neighboring
points $t_1 < t_2$ of different groups $J_m, J_{m+1} \subset J$.
Then groups $\Xi_*(\J), \Xi_*(\J')$ are obviously identified. Let
$A!$ be the multiindex, obtained from $A$ by replacing two numbers
$a_m, a_{m+1}$ by one number $a_m+a_{m+1}-1$, and $J!$ the $A!$-configuration
obtained from $J$ by contracting the segment $[t_1,t_2]$ into
a point $\tau$, which
will thus belong to the group of cardinality $a_m+a_{m+1}-1$.
The index
$\langle V|J! \rangle$
is an element of the group $\Xi_*(J!)$,
i.e. a linear combination of collections of $\#A!$ two-connected
graphs (with ordered edges). Its boundary
$\partial
\langle V|J! \rangle$
is a linear combination of similar collections of graphs, exactly one
of which is not two-connected. Denote by
$\partial_\tau
\langle V|J! \rangle$
the part of this linear combination, spanned by all such collections,
that removing from any of them the point $\tau$ we split the corresponding
graph $\Gamma_m$ with $a_m+ a_{m+1}-1$ vertices into two disjoint graphs, whose
nodes are the points of $J_m \sm t_1$ and $J_{m+1} \sm t_2$. Any such graph
$\Gamma_m$ is the union of two graphs with $a_m$ and $a_{m+1}$ nodes and
common node $\tau$, thus  the space of all such linear
combinations also can be identified with any of groups
$\Xi_*(\J), \Xi_*(\J').$ The promised homological condition
consists in the fact that
the chain in $\sigma_i$, Alexander dual to $V$, actually is a cycle,
and in particular its algebraic boundary close to the common boundary
of $\J$-, $\J'$- and $\J!$-blocks is equal to zero. This condition is
as follows:
\begin{equation}
\label{stu}
\langle V|\J \rangle -
\langle V|\J' \rangle =
\pm \partial_\tau \langle V|\J! \rangle,
\end{equation}
where the coefficient $\pm$ depends on the choice of local
orientations of all participating strata.

The system of equations (\ref{stu}), corresponding to all
points $\tau$ of the configuration $\J!$,
is strong enough to determine the index
$\langle V|\J! \rangle$ if we know all indices
$\langle V|\J \rangle$ for all configurations $\J$ of smaller
stickiness. Indeed, any knot invariant $V$
of order $i$ is determined (up to terms of lower orders) by
its indices at chord diagrams, i.e., the configurations of
stickiness $0$.

This equality allows us also to give the following
characterization of the reversed filtration. Let $\J$, $\J'$ be
two classes of equivalent $A$-configurations with the same $A$.
An arbitrary correspondence between their groups $J_m, J'_{m'}$ of
equal cardinalities allows us to identify the groups
$\Xi_*(\J)$ and $\Xi_*(\J')$.
This identification is not unique if
some of numbers $a_m$ of the multiindex $A$ coincide. Moreover,
in the last case we can define a group of automorphisms, acting
on $\Xi_*(\J)$ and generated by all possible permutations of
groups of the same cardinalities $a_m$, preserving the order of points inside
any group, cf. Corollary 1 above.
\medskip

{\sc Proposition 5.} {\it For any knot invariant $V$ of order $i$,
the following conditions are
equivalent:

1) $V$ has reversed filtration $\le \al$;

2) for any multiindex $A$ with $|A|-\# A=i,$ $i-\# A = \al,$
all the indices $\langle V|\J\rangle \in \Xi_*(\J)$, defined
by $V$ in all $\J$-blocks over all $A$-configurations $\J,$
are invariant under all the identifications and automorphisms
described in the previous paragraph.}
\medskip

{\it Proof.} (1) $\Rightarrow$ (2).
If $\bar V$ is the class in $\bar H_{\omega-1}(\sigma_i)$, representing
$V$, then for any $A$ as above
its restriction on $\Phi_\al \sm \Phi_{\al-1}$
defines a class in $\bar H_{\omega-1}(S(A))$.
The indices $\langle V|\J \rangle$ for all classes $\J$ of equivalent
$A$-configurations
are defined by
restrictions of this class to all corresponding
$\J$-blocks. In particular, these indices
for neighboring classes of
$A$-configurations (i.e. for those obtained one from another by a
change of orders of only two points from different groups) should coincide:
otherwise our cycle in $S(A)$ would have nontrivial boundary
at the border between these two blocks.

The implication (2) $\Rightarrow$ (1) means
that there are no non-trivial elements of the group
$\bar H_{\omega-1}(\sigma_i \sm \sigma_{i-1})$, represented by
cycles with support in $(\sigma_i \sm \sigma_{i-1}) \sm \Phi_\al$;
this  follows from the fact
that any invariant $V$ of order $i$ is determined up to invariants
of lower orders by its indices at all $i$-chord diagrams
(i.e., at all $\J$-blocks with $A=(2, \ldots, 2)$). \quad $\Box$

\subsection{The symbol of a cohomology class of a finite order}

Let again $n$ be an arbitrary natural number greater than 2,
and  $K$ the space of non-compact knots $\R^1 \to \R^n$. Let
$V \in H^a(K \sm \Sigma)$ be any cohomology class of order $i$, i.e.,
a class, Alexander dual to an $(\omega-a-1)$-dimensional cycle $\bar V \in
\sigma_i$.  Let $\al$ be the reversed filtration of this class, i.e. we can
choose this cycle $\bar V$ in such a way that it lies in
$\Phi_\al \cup \sigma_{i-1}$.
Then for an arbitrary multiindex $A$ with $|A|-\#A=i$ and $\#A=i-\al$,
the {\it symbol} $s(A,V)$ is defined as follows.

Let us denote by $\Xi(A)$ the local system of groups on the configuration
space $C(A)$, whose fiber over an $A$-set $\U$ is identified with $\Xi_*(\U)$.
This local system is not constant if some of numbers
$a_m$ coincide: a path in $C(A)$, permuting some groups of the
same cardinality, acts nontrivially on the fiber.

Further, let $\tilde \Xi(A)$ be a similar local system with a slightly more
complicated monodromy action. Namely, any path in $C(A)$, permuting exactly
two groups of points, acts exactly as in the system $\Xi(A)$ (i.e., permutes
the corresponding factors $\bar H_*(\Lambda(J_m))$ ) if the cardinality
of any of these groups is odd, and additionally multiplies such an operator
by $-1$ if these cardinalities are even.

The symbol $s(A,V)$, which we are going to define, is an element of the group
$H^a(C(A),\Xi(A))$ if $n$ is odd and of
$H^a(C(A),\tilde \Xi(A))$ if $n$ is  even.

The construction repeats that for the index of a knot invariant given in the
previous subsection. Namely, we consider the class of the realizing
$V$ dual cycle $\bar V$ in the group $\bar H_{\omega-a-1}(\Phi_\al \sm
\Phi_{\al-1})$, then, using Corollary 2 of Theorem 3, Thom isomorphism
for the fiber bundle mentioned in this Corollary,
and Poincar\'e duality, we reduce this group to the cohomology group
of $C(A)$ with coefficients in a local system associated to the homology
bundle of the fiber bundle $
\Phi_\al \sm
\Phi_{\al-1}\to C(A)$.

For instance, if $n=3$ and $a=0$, i.e. $V$ is just a knot invariant, then
this symbol coincides with the totality of indices of $V$, corresponding
to all classes of $A$-configurations. The invariance of these indices,
mentioned in Proposition 5, follows from the fact, that
the 0-dimensional cocycle
of a local system is just a global section of this system and thus
is invariant under its monodromy action.

\subsection{Multiplication}

\subsubsection{Multiplication of knot invariants in $\R^3$.}
We have two filtrations in the space of finite-order invariants of knots
in $\R^3$: the order $i$ and, for any fixed $i,$
the degree $\al$ with respect to the reversed
auxiliary filtration in $\sigma_{i} \sm \sigma_{i-1}.$
It is well-known that the order is multiplicative. In fact, the same
is true also for the second number. \medskip

{\sc Theorem 4.} {\it Let $V',V''$ be two knot invariants of orders
$i', i''$. Then
for any multiindex $A=(a_1, \ldots, a_{\# A})$
of complexity $i=i'+i''$
and any equivalence class $\J$ of $A$-configurations,
the index $\langle V'\cdot V''|\J\rangle \in
\Xi_*(\J)$ is equal
to \begin{equation} \label{tens}
\sum_{J' \cup J''=J} \langle V'|J'\rangle \otimes \langle V''|J''\rangle,
\end{equation} where summation is taken over all decompositions of
any $A$-configuration $J$ representing $\J$ into
disjoint unions of an $A'$-configuration $J'$ and $A''$-configuration $J''$ of
complexities $i'$ and $i''$, such that any group of $J$
belongs to either $J'$ or $J'$.
} \medskip

{\it Proof.} For
$A=(2^{\times i})$ this is just the Kontsevich's
decomposition formula for the index of $V' \cdot V''$ on an
$i$-chord diagram:
for any $(2, \ldots, 2)$-configuration $\tau$,
\begin{equation}
\label{konts}
\langle(V' \cdot V''),(\tau)\rangle =
\sum_{\tau' \cup \tau''= \tau} V'(\tau')V''(\tau''),
\end{equation}
summation over all $\binom{i}{i'}$ ordered decompositions of the
$i$-chord diagram $\tau$ into an $i'$-chord diagram $\tau'$
and an $i''$-chord diagram $\tau''$; see e.g. \cite{BN2}.

For an arbitrary multiindex $A$ of the same complexity $i$
the formula (\ref{tens}) follows by induction over the stickiness, while the
induction step is the formula (\ref{stu}). \quad $\Box$
\medskip

{\sc Corollary 4.} {\it The bi-order $(i,\al)$ is multiplicative:
if $V'$ and $V''$ are two invariants of bi-orders $(i',\al')$ and
$(i'',\al'')$ respectively, then $V' \cdot V''$ is an invariant of bi-order
$(i'+i'',\al'+\al'').$} \medskip

Indeed, in this case the formula (\ref{tens}) gives one and the
same answer for all $A$-configurations $J$ of complexity $i'+i''$
and stickiness $\alpha'+\alpha''$, thus our Corollary follows
from Proposition 5. \quad $\Box$

%\medskip
%{\sc Corollary 2.} {\it The group $\Phi_\alpha$ of invariants
%of bi-order $(i,\alpha)$ is generated (modulo the group of invariants
%of orders $<i$) by all products of $i-\al$ invariants of positive  orders
%$i_1, \ldots, i_{i-\al}$, $i_1 + \ldots + i_{i-\al}= i$.} \medskip

%In particular, Proposition 8$'$ gives us an upper estimate of the number of
%order $i$ independent multiplicative generators of the algebra of invariants
%by the dimension of the subgroup of $\Z_{i+1}$-stable elements of
%$\bar H_{2i-2}(\Delta^2(i+1))$.

%{\it Proof.}

\subsubsection{Multiplication conjecture for higher cohomology
of spa\-ces of knots}
For any miltiindices $A', A''$ and $A \equiv A' \cup A'',$ there
are natural operations
\begin{equation}
H^*(C(A'),\Xi(A')) \otimes H^*(C(A''),\Xi(A'')) \to H^*(C(A),\Xi(A)),
\label{conj1}
\end{equation}
\begin{equation}
H^*(C(A'),\tilde \Xi(A')) \otimes H^*(C(A''),\tilde \Xi(A'')) \to
H^*(C(A),\tilde \Xi(A))
\label{conj2}
\end{equation}
such that for any cohomology classes $V', V'' \in H^*(\K \sm \Sigma)$
of bi-orders $(i',\al')$ and $(i'',\al'')$ respectively,
the class $V' \cdot V''$ has bi-order $(i'+i'', \al'+\al''),$ and
its principal symbol $s(A, V' \cdot V'')$ can be expressed through
the symbols of $V'$ and $V''$ by a formula similar to
(\ref{tens}), in which the operation of the tensor multiplication
is replaced by the operation (\ref{conj1}) if $n$ is odd
and by operation (\ref{conj2}) if $n$ is even, cf.
\cite{Fu}, \S \ 8.

\subsection{Configurations, containing groups of 2 points,
do not contribute to the cohomology of lowest possible dimension.}

Consider the {\it reversed spectral sequence} $\E^r_{p,q}$,
calculating the group $\bar H_*(\sigma_i \sm \sigma_{i-1})$ and
generated by the reversed filtration $\{\Phi_\al\}.$
By definition, the term $\E^1_{\al,q}$ of this spectral sequence
is isomorphic to $\bar H_{\al+q}(\Phi_\al \sm \Phi_{\al-1})$ and thus splits
into the direct sum of similar groups $\bar H_{\al+q}(S(A))$ over all indices
$A$ of the form (\ref{a}) with $\# A =i-\al.$

In particular, by Corollary 3 it is trivial if $\al+q> \omega-1-(n-3)i.$
\medskip

For instance, for any $i$ the term $\Phi_0$ of our filtration
consists of unique stratum $S(A)$ with $A=(2, \ldots, 2).$
\medskip

{\sc Proposition 6.} {\it For any $i$, two top terms
$\E^1_{0,\omega-1-(n-3)i}$ and $\E^1_{0,\omega-2-(n-3)i}$ of the column
$\E^1_{0,*} \equiv \E^1_{0,*}(i)$ of the reversed auxiliary spectral
sequence, calculating $\bar H_*(\sigma_i \sm \sigma_{i-1})$, are trivial.}
\medskip

This proposition is a special case of the following one.
\medskip

{\sc Proposition 7.} {\it For any multiindex $A$ of the form (\ref{a}),
such that at least one (respectively, at least two)
of numbers $a_i$ are equal to 2, the group
$\bar H_{\omega-1-(n-3)i}(S(A))$ (respectively, both groups
$\bar H_{\omega-1-(n-3)i}(S(A)),$ $\bar H_{\omega-2-(n-3)i}(S(A))$)
is trivial.}
\medskip

{\it Example.} Let be $n=3$, then, calculating
the group $\bar H_{\omega-1}(\sigma_i)$
of knot invariants of order $i$, we can take into account only
strata $S(A)$ with at most one character 2 in the multiindex $A.$
Moreover, such strata with exactly one 2 in $A$ can provide only
relations in this group, and not generators.
\medskip

{\sc Lemma 3.} {\it For any natural $k$,
$\bar H_*(B(\Psi,k)) \simeq 0$ and $\bar H_*(B(\Psi,k),\pm \Z) \simeq 0.$}
\medskip

(For the definition of the local system $\pm Z$, see \S \ 5.1 above.)
\medskip

{\it Proof of Lemma 3.} We use the decomposition of $B(\Psi,k)$ similar to
the cell decomposition of the space $B(\R^2,k)$ used in
\cite{Fu}. Namely, to any decomposition $k=k_1+ \cdots+k_m$ of
the number $k$ into natural numbers $k_i$ we assign the set of all
$k$-subsets of $\Psi$, consisting of points
$((t_1\le t'_1), \ldots, (t_k  \le t'_k)),$ such that the smallest value
of numbers $t_l$ appears in exactly $k_1$ pairs $(t_l \le t'_l),$
the next small value in $k_2$ pairs etc. Denote this set by
$e(k_1, \ldots, k_m)$. Filter the space $B(\Psi,k)$, assigning to
the $l$-th term of the filtration the union of all such sets
of dimension $\le l.$ It is easy to see that
any such set is diffeomorphic to the direct product of a
closed octant in $\R^m$ and an open $k$-dimensional cell,
thus its Borel--Moore homology group is trivial, as well as
the similar group with coefficients in (the restriction on this set of) the
local system $\pm \Z.$
Therefore the spectral sequences, calculating both groups
$\bar H_*(B(\Psi,k)), \bar H_*(B(\Psi,k),\pm \Z)$ and generated
by our filtration, vanish in the
first term. \quad $\Box$
\medskip

{\it Proof of Proposition 6.} If $A=(2^{\times i})$, then
$C(A)$ is an open subset in the configuration space $B(\Psi,i)$, whose
complement $B(\Psi,i) \sm C(A)$ has codimension 3. In particular,
by Lemma 3 $\bar
H_j(C(A))=0$ for $j=2i$ or $2i-1.$ It is easy to calculate that
for odd $n$ the fiber bundle
$S(A) \to C(A)$ is orientable, and for even $n$ it changes its orientation
together with the local system $\pm \Z$. Thus we have the
Thom isomorphism
$\bar H_{*+\omega-(n-2)i-1}(S(A)) \simeq \bar H_{*}(C(A))$
(for odd $n$) and $\simeq \bar H_{*}(C(A),\pm \Z)$ (for even $n$).
\quad $\Box$ \medskip

{\sc Lemma 4.} {\it For any finite subset $\theta \subset \Psi$, both
groups $\bar H_l(B(\Psi \sm \theta, k)), $
$\bar H_l(B(\Psi \sm \theta,k),\pm \Z)$ are trivial for $l \ne k,$
and for $l=k$ they are free Abelian of rank equal to the number of functions
$\chi: \theta \to \Z_+$ such that $\sum_{z \in \theta} \chi(z)=k.$}
\medskip

{\it Proof.} Consider any direction in the plane
$\R^2 \supset \Psi$, transversal to
the diagonal $\partial \Psi = \{t=t'\}$ and such that none two
points of $\theta$ are connected by a vector of this direction.
Let $\pi: \Psi \to \partial \Psi$ be the projection along this
direction. To any point $Z = (z_1, \ldots, z_k) \subset \Psi \sm \theta$
of the space $B(\Psi \sm \theta, k)$ the following data are assigned:

(1) the topological type of the configuration in $\partial \Psi,$
formed by the points $\pi(z_j)$ (counted with multiplicities)
and the set $\pi(\theta);$

(2) for any point $w \in \theta,$ the number of points $z_j \in Z$ such that
$\pi(z_j) = \pi(w)$ and $z_j$ is separated by $w$ from $\partial \Psi$
in the line $\pi^{-1}(\pi(w))$.

For any such collection of data,
the subset in $B(\Psi \sm \theta,k)$,
formed by configurations $Z$ with these data, is homeomorphic to
$\R^l_+ \times \R^s,$ where $l$ is the number of lines
$\pi^{-1}(\cdot)$ containing points of $Z$, not separated from $\partial \Psi$
by the points of $\theta,$ and $l+s$ is equal to the dimension of this
subset, i.e. to $k+$(the number of geometrically distinct points of
$\pi(Z)$ for any $Z$ from this subset). In particular, the
Borel--Moore homology group of such a subset is trivial if $l>0,$
and the set of configurations with $l=0$ consists of
several $k$-dimensional cells. \quad $\Box$ \medskip

{\it Proof of Proposition 7.} If exactly one of numbers $a_j$ is equal
to 2, then $S(A)$ is a connected manifold with non-empty boundary,
in particular its Borel--Moore homology group
of top dimension is trivial.

Now suppose that there are at least two such numbers.
Let $A_{>2}$ be the same multiindex $A,$
from which all numbers $a_j$ equal to 2 are removed. Consider the projection
\begin{equation}
C(A) \to C(A_{>2}),
\label{b2}
\end{equation}
erasing from any $A$-set all its groups of cardinality 2. For any point
$\nu \in C(A_{>2})$ denote by $\theta(\nu)$ the set of all pairs
$(x,y)$ lying in some of groups of $\nu$.

By the Leray spectral sequence of the composite fiber bundle
$S(A) \to C(A) \to C(A_{>2}),$ we need only to prove the following lemma.
\medskip

{\sc Lemma 5.} {\it If the number $r$ of twos in $A$ is greater than
1, then for any $\nu \in C(A_{>2})$ the group
$\bar H_*$ of the fiber of the projection (\ref{b2})
over the point $\nu$ is trivial in
dimensions $2r$ and $2r-1$.}
\medskip

{\it Proof.} Any such fiber consists of all configurations
$Z \in B(\Psi,r),$  not containing the points of the
finite subset $\theta(\nu) \subset \Psi$ and satisfying some additional
restrictions, which forbid certain subvariety of codimension $\ge 3$
in $B(\Psi,r).$ Thus our lemma follows from Lemma 4. \quad $\Box$
\medskip

Recall that the greatest possible value of the reversed filtration
of a cohomology class (in particular, of an knot invariant in $\R^3$)
of order $i$ is equal to $i-1.$
\medskip

{\sc Proposition 8.} {\it For any $i$,
the number of linearly independent knot invariants of
bi-order $(i,i-1)$
(modulo invariants of lower bi-orders)
is estimated from above by the number
$(i-1)! \equiv \dim \bar H_{2i-2}(\Delta^2(i+1)).$}
\medskip

This follows from Theorem 2 and the fact that unique multiindex $A$
of the form (\ref{a}) with given $i$ and stickiness $i-1$
consists of one number $(i+1).$
\quad $\Box$ \medskip

This estimate is not realistic. Indeed, the group $\Z_{i+1}$
of cyclic permutations of vertices acts naturally on the complex
$\Delta(i+1)$, hence also on the group $\bar H_*(\Delta^2(i+1))$.
It is easy to see, that an element of this group, corresponding to a knot
invariant of bi-order $(i,i-1)$, should be invariant under this action,
and we obtain the following improvement of Proposition 8.
\medskip

{\sc Proposition 8$'$.} {\it
The number of linearly independent knot invariants of
bi-order $(i,i-1)$ (modulo invariants of lower bi-orders)
does not exceed the dimension of the subgroup
in $\bar H_{2i-2}(\Delta^2(i+1))$, consisting of $\Z_{i+1}$-invariant
elements.} \quad $\Box$

The exact formula for this dimension was found in \cite{BBLSW}, see
Corollary 4.7 there.

\section{Order 3 cohomology of spaces of non-compact knots}

Here we calculate the column $E^\infty_{3,*}$ of the main spectral sequence,
converging to the finite-order cohomology of the space of non-compact knots
$\R^1 \to \R^n$. In fact, it is sufficient to calculate the column
$E^1_{3,*}$ of its initial term $E^1$, i.e. the group
$\bar H_*(\sigma_3 \sm \sigma_2)$
of the corresponding discriminant variety.

This group is described in the following statement.
\medskip

{\sc Theorem 5.} {\it For any $n \ge 3$, all groups
$\bar H_j(\sigma_3 \sm \sigma_2)$
are trivial, except for such groups with
$j=\omega-1-3(n-3)$ and $\omega-2-3(n-3),$
which are isomorphic to $\Z$.}
\medskip

This statement for $n=3$ is not new: the group
$\bar H_{\omega-1}(\sigma_3 \sm \sigma_2) \simeq \Z$
(of knot invariants of order 3) was calculated
in \cite{V1}, and all other groups
$\bar H_j(\sigma_3 \sm \sigma_2)$
in a non-published work of D.~M.~Teiblum and V.~E.~Turchin,
see \S \ 4.4.2 above.
Their calculation is based on the cellular decomposition of
$\sigma_3 \sm \sigma_2,$ constructed in \cite{V1}, and is non-trivial
even for a computer.
\medskip

{\sc Corollary 5.} {\it The column $E^\infty_{3,*}$ of the main
spectral sequence coincides with $E^1_{3,*}$, namely, it consists of
exactly two nontrivial terms
$E^\infty_{3,\omega-1-3(n-2)} \simeq \Z$ and
$E^\infty_{3,\omega-2-3(n-2)} \simeq \Z.$}
\medskip

{\it Proof of Corollary 5.} The fact that all differentials
$d^r : E^r_{3,*} \to E^r_{3-r,*+r-1},$  $r \ge 1,$ are trivial,
follows immediately from the construction of columns
$E^1_{2,*}$ and $E^1_{1,*},$ see Examples 3 and 4 in \S \ 5.1.

On the other hand, the inequality (\ref{ineq}) implies that also there are no
non-trivial differentials $d^r$ acting {\it into} the cell
$E^r_{3,\omega-1-3(n-2)}$; if $n >3$ then the same is true also for the
cell $E^r_{3,\omega-2-3(n-2)}$.  Finally, if $n=3,$ then the similar
triviality of the homomorphisms
$d^r:E^r_{3+r, \omega-4-r} \to E^r_{3,\omega-5}$  follows from the
Kontsevich's realization theorem. \quad $\Box$

\subsection{Term $E^1$ of the reversed spectral sequence.}

The term $\sigma_3 \sm \sigma_2$ of the main filtration consists of
three terms of the reversed auxiliary filtration:
$\Phi_0 \subset \Phi_1 \subset \Phi_2 \equiv \sigma_3 \sm \sigma_2.$

The sets $\Phi_0,$ $\Phi_1 \sm \Phi_0$ and $\Phi_2 \sm \Phi_1$
consist of one component $S(A)$ each, with $A$ equal to $(2,2,2),$
$(3,2)$ and $(4),$ respectively.
\medskip

{\sc Lemma 6}. {\it For any $n$,

1. $\bar H_j(\Phi_0)=0$ for all $j \ne \omega-3n+6,$
and
$\bar H_{\omega-3n+6}(\Phi_0)=\Z$.

2. $\bar H_j(\Phi_1 \sm \Phi_0)=0$ for all $j \ne \omega-3n+7,$
and
$\bar H_{\omega-3n+7}(\Phi_1 \sm \Phi_0)=\Z^3$.

3. $\bar H_j(\Phi_2 \sm \Phi_1)=0$ for all $j \ne \omega-3n+8,$
and $\bar H_{\omega-3n+8}(\Phi_2 \sm \Phi_1)=\Z^2$.}
\medskip

{\it Proof.} 1. The space $C(2,2,2)$ of all $(2,2,2)$-sets
is an open subset in the configuration space $B(\Psi,3)$, and their
difference ${\bf T} \equiv B(\Psi,3)\sm C(2,2,2)$ is the set of all triples
of the form $(x,y),(x,z),(y,z) \subset \Psi$ with $x<y<z.$
This set is obviously a closed subset in $B(\Psi,3)$
diffeomorphic to a 3-dimensional cell. Thus by Lemma 3
and the exact sequence of the pair $(B(\Psi,3),{\bf T})$,
both groups $\bar H_l(C(2,2,2))$ and $\bar H_l(C(2,2,2),\pm \Z)$
are trivial for any $l$ other than $4$
and are isomorphic to $\Z$ in dimension $4$.
Fibres of the bundle $S(2,2,2) \to C(2,2,2)$ are products
of the open triangle $\Lambda(2,2,2)$ and an affine space
of dimension $\omega-3n$.
It is easy to calculate that for odd $n$ this bundle is orientable, and for
even $n$ its orientation bundle coincides with $\pm \Z$,
thus statement 1 follows from the Thom isomorphism.

\begin{figure}
\unitlength=1.00mm
\special{em:linewidth 0.4pt}
\linethickness{0.4pt}
\begin{picture}(53.00,52.00)
\thicklines
\put(5.00,1.00){\line(1,1){45.00}}
\put(10.00,6.00){\vector(1,0){40.00}}
\put(10.00,6.00){\vector(0,1){40.00}}
\thinlines
\put(7.00,46.00){\makebox(0,0)[cc]{$t'$}}
\put(53.00,3.00){\makebox(0,0)[cc]{$t$}}
\put(20.50,5.50){\rule{1.00\unitlength}{1.00\unitlength}}
\put(28.50,5.50){\rule{1.00\unitlength}{1.00\unitlength}}
\put(29.00,2.00){\makebox(0,0)[cc]{$y$}}
\put(21.00,2.00){\makebox(0,0)[cc]{$x$}}
\put(7.00,23.00){\makebox(0,0)[cc]{$y$}}
\put(9.50,24.50){\rule{1.00\unitlength}{1.00\unitlength}}
\put(9.50,35.50){\rule{1.00\unitlength}{1.00\unitlength}}
\put(7.00,34.00){\makebox(0,0)[cc]{$z$}}
\put(21.00,25.00){\circle{2.00}}
\put(21.00,36.00){\circle{2.00}}
\put(29.00,36.00){\circle{2.00}}
\put(29.00,46.00){\vector(0,-1){9.00}}
\put(12.00,45.00){\vector(1,-1){8.00}}
\put(5.00,25.00){\vector(1,0){15.00}}
\put(5.00,25.00){\line(-1,0){3.00}}
\put(12.00,45.00){\line(-1,1){7.00}}
\put(29.00,52.00){\line(0,-1){6.00}}
\end{picture}
\caption{Basic cycles in the space $C(3,2)$}
\label{aaa5}
\end{figure}
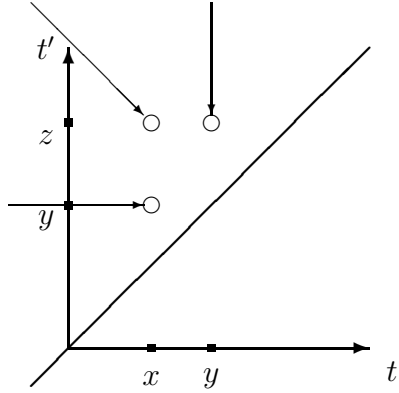

2. The space $C(3,2)$ can be considered as the space of a fiber bundle,
whose base is the configuration space $B(\R^1,3)$, and the
fiber over the triple $(x<y<z) \subset \R^1$ is the half-plane
$\Psi$ with three interior points $(x,y), (x,z)$ and $(y,z)$ removed.
The Borel--Moore homology group of this punctured half-plane is
concentrated in dimension 1 and is generated by three rays connecting
these three removed points to infinity, say, by rays shown in
Fig.~\ref{aaa5}.
The complex $\Lambda(3,2)$ consists of unique
3-dimensional open face. The bundle $S(3,2) \to C(3,2)$
is orientable, thus statement 2 also follows from the Thom
isomorphism.

3. The space $C(4) \sim B(\R^1,4)$ is obviously diffeomorphic
to a 4-cell, thus $S((4))$ is homeomorphic to the direct product
$B(\R^1,4) \times \R^{\omega-3n} \times \Delta^2(n)$, and  statement 3
follows immediately from Theorem 2. For the realization of generators
of the group $\bar H_4(\Delta^2(4)) \simeq
\bar H_{\omega-3n+8}(\Phi_2 \sm \Phi_1)$
see Example 2 in \S \ 2. \quad $\Box$
\medskip

{\sc Remark.} To specify the last isomorphism, we need
to choose an orientation of
$B(\R^1,4)$ and a (co)orientation of the fiber $\R^{\omega-3n}.$
We do it as follows.

The standard orientation of $B(\R^1,4)$   is defined by
local coordinate systems, whose coordinate functions are
coordinates of four points taken in their increasing order.

If $\xi_1, \ldots, \xi_n$ are coordinates in $\R^n$, so that any
non-compact knot is given
parametrically by $n$ real functions $\xi_k(t)$, $k=1, \ldots, n$, then the
fiber $\R^{\omega-3n}$ over the configuration $(r<s<t<u)$ is distinguished by
the equations $\xi_1(r)=\xi_1(s)=\xi_1(t)=\xi_1(u)$, \ldots,
$\xi_n(r)=\xi_n(s)=\xi_n(t)=\xi_n(u).$ Then the transversal orientation of this
fiber in the space $K \sim \R^\omega$ will be specified by the skew form
\begin{eqnarray} \label{frame}
\quad d(\xi_1(s)-\xi_1(r)) \wedge d(\xi_2(s)-\xi_2(r)) \wedge \ldots \wedge
d(\xi_n(s)-\xi_n(r)) \wedge \nonumber \\
\wedge \ d(\xi_1(t)-\xi_1(s)) \wedge d(\xi_2(t)-\xi_2(s)) \wedge \ldots \wedge
d(\xi_n(t)-\xi_n(s)) \wedge \nonumber \\
\wedge \ d(\xi_1(u)-\xi_1(t)) \wedge d(\xi_2(u)-\xi_2(t)) \wedge \ldots \wedge
d(\xi_n(u)-\xi_n(t)). \
\end{eqnarray}
\smallskip

Finally, the calculation of groups $\bar H_j(\sigma_3 \sm \sigma_2)$
in dimensions $j=\omega-3n+8,$ $\omega-3n+7$ and $\omega-3n+6$ is reduced
to the calculation of a certain complex of the form
$0 \to \Z^2 \to \Z^3 \to \Z \to 0,$ whose boundary operator
is the (horizontal) differential $d_1$ of our reversed spectral sequence.

\subsection{The differential
$d^1: \bar H_{\omega-3n+8}(\Phi_2 \sm \Phi_1) \to
\bar H_{\omega-3n+7}(\Phi_1 \sm \Phi_0)$.}
First we calculate this operator for the generator
of the group
$\bar H_{\omega -3n+8}(\Phi_2 \sm \Phi_1)$, corresponding to the
generator

\begin{equation}
\label{first}
\unitlength=1.00mm
\special{em:linewidth 0.4pt}
\linethickness{0.4pt}
\begin{picture}(28.00,4.00)
\put(2.00,3.00){\makebox(0,0)[cc]{1}}
\put(2.00,-1.00){\makebox(0,0)[cc]{2}}
\put(12.00,-1.00){\makebox(0,0)[cc]{3}}
\put(12.00,3.00){\makebox(0,0)[cc]{4}}
\put(28.00,3.00){\makebox(0,0)[cc]{4}}
\put(28.00,-1.00){\makebox(0,0)[cc]{3}}
\put(18.00,-1.00){\makebox(0,0)[cc]{2}}
\put(18.00,3.00){\makebox(0,0)[cc]{1}}
\put(15.00,1.00){\makebox(0,0)[cc]{{\large $-$}}}
\put(4.00,-2.00){\line(1,0){6.00}}
\put(10.00,-2.00){\line(0,1){6.00}}
\put(10.00,4.00){\line(-1,0){6.00}}
\put(4.00,4.00){\line(0,-1){6.00}}
\put(4.00,-2.00){\line(1,1){6.00}}
\put(20.00,4.00){\line(0,-1){6.00}}
\put(20.00,-2.00){\line(1,0){6.00}}
\put(26.00,-2.00){\line(0,1){6.00}}
\put(26.00,4.00){\line(-1,0){6.00}}
\put(20.00,4.00){\line(1,-1){6.00}}
\end{picture}
\end{equation}

\noindent
of $\bar H_4(\Delta^2(4))$,
where numbers of vertices correspond to the order
of four points in $\R^1.$
The boundary of any of these two graphs is the sum of 5 graphs
with 4 edges. One of them (obtained by removing the diagonal edge)
appears in both sums
and vanishes. Remaining boundary graphs (and the signs with which
they participate in the algebraic boundary of the cycle (\ref{first})) are:

\begin{equation}
\label{aaa3}
\unitlength=1.00mm
\special{em:linewidth 0.4pt}
\linethickness{0.4pt}
\begin{picture}(136.00,12.00)
\put(121.00,9.00){\oval(26.00,6.00)[b]}
\put(112.50,9.00){\oval(9.00,4.00)[t]}
\put(125.50,9.00){\oval(17.00,6.00)[t]}
\put(121.00,9.00){\oval(8.00,4.00)[b]}
\put(86.00,9.00){\oval(26.00,6.00)[b]}
\put(77.50,9.00){\oval(9.00,4.00)[t]}
\put(90.50,9.00){\oval(17.00,6.00)[t]}
\put(94.50,9.00){\oval(9.00,4.00)[b]}
\put(55.50,9.00){\oval(17.00,6.00)[b]}
\put(51.00,9.00){\oval(8.00,4.00)[t]}
\put(59.50,9.00){\oval(9.00,4.00)[t]}
\put(20.50,9.00){\oval(17.00,6.00)[b]}
\put(16.00,9.00){\oval(8.00,4.00)[t]}
\put(24.50,9.00){\oval(9.00,4.00)[t]}
\put(42.50,9.00){\oval(9.00,4.00)[b]}
\put(16.00,9.00){\oval(26.00,6.00)[t]}
\put(121.00,2.00){\makebox(0,0)[cc]{D($-$)}}
\put(86.00,2.00){\makebox(0,0)[cc]{C($-$)}}
\put(51.00,2.00){\makebox(0,0)[cc]{B(+)}}
\put(16.00,2.00){\makebox(0,0)[cc]{A($-$)}}
\put(1.00,9.00){\line(1,0){30.00}}
\put(36.00,9.00){\line(1,0){30.00}}
\put(71.00,9.00){\line(1,0){30.00}}
\put(106.00,9.00){\line(1,0){30.00}}
\put(3.00,9.00){\circle*{1.33}}
\put(12.00,9.00){\circle*{1.33}}
\put(20.00,9.00){\circle*{1.33}}
\put(29.00,9.00){\circle*{1.33}}
\put(38.00,9.00){\circle*{1.33}}
\put(47.00,9.00){\circle*{1.33}}
\put(55.00,9.00){\circle*{1.33}}
\put(64.00,9.00){\circle*{1.33}}
\put(73.00,9.00){\circle*{1.33}}
\put(82.00,9.00){\circle*{1.33}}
\put(90.00,9.00){\circle*{1.33}}
\put(99.00,9.00){\circle*{1.33}}
\put(108.00,9.00){\circle*{1.33}}
\put(117.00,9.00){\circle*{1.33}}
\put(125.00,9.00){\circle*{1.33}}
\put(134.00,9.00){\circle*{1.33}}
\end{picture}
\end{equation}

\noindent
for the first graph in (\ref{first}), and

\begin{equation}
\label{aaa4}
\unitlength=1.00mm
\special{em:linewidth 0.4pt}
\linethickness{0.4pt}
\begin{picture}(135.00,13.00)
\put(15.00,10.00){\oval(26.00,6.00)[b]}
\put(10.50,10.00){\oval(17.00,6.00)[t]}
\put(23.50,10.00){\oval(9.00,4.00)[t]}
\put(85.00,10.00){\oval(26.00,6.00)[b]}
\put(80.50,10.00){\oval(17.00,6.00)[t]}
\put(93.50,10.00){\oval(9.00,4.00)[t]}
\put(76.50,10.00){\oval(9.00,4.00)[b]}
\put(15.00,10.00){\oval(8.00,4.00)[b]}
\put(45.50,10.00){\oval(17.00,6.00)[b]}
\put(41.50,10.00){\oval(9.00,4.00)[t]}
\put(50.00,10.00){\oval(8.00,4.00)[t]}
\put(115.50,10.00){\oval(17.00,6.00)[b]}
\put(111.50,10.00){\oval(9.00,4.00)[t]}
\put(120.00,10.00){\oval(8.00,4.00)[t]}
\put(120.00,10.00){\oval(26.00,6.00)[t]}
\put(58.50,10.00){\oval(9.00,4.00)[b]}
\put(15.00,3.00){\makebox(0,0)[cc]{E(+)}}
\put(50.00,3.00){\makebox(0,0)[cc]{F(+)}}
\put(85.00,3.00){\makebox(0,0)[cc]{G($-$)}}
\put(120.00,3.00){\makebox(0,0)[cc]{H(+)}}
\put(0.00,10.00){\line(1,0){30.00}}
\put(35.00,10.00){\line(1,0){30.00}}
\put(70.00,10.00){\line(1,0){30.00}}
\put(105.00,10.00){\line(1,0){30.00}}
\put(2.00,10.00){\circle*{1.33}}
\put(11.00,10.00){\circle*{1.33}}
\put(19.00,10.00){\circle*{1.33}}
\put(28.00,10.00){\circle*{1.33}}
\put(37.00,10.00){\circle*{1.33}}
\put(46.00,10.00){\circle*{1.33}}
\put(54.00,10.00){\circle*{1.33}}
\put(63.00,10.00){\circle*{1.33}}
\put(72.00,10.00){\circle*{1.33}}
\put(81.00,10.00){\circle*{1.33}}
\put(89.00,10.00){\circle*{1.33}}
\put(98.00,10.00){\circle*{1.33}}
\put(107.00,10.00){\circle*{1.33}}
\put(116.00,10.00){\circle*{1.33}}
\put(124.00,10.00){\circle*{1.33}}
\put(133.00,10.00){\circle*{1.33}}
\end{picture}
\end{equation}

\noindent
for the second. The term $\Phi_1 \sm \Phi_0$ (in which the boundary lies)
can be considered as a fiber bundle (we will call it the {\it former bundle}),
whose base is the set of triples $\{x<y<z\} \subset \R^1,$
and the fiber over such a point is the space of the {\it latter fiber bundle},
whose base is the set of all pairs $\{u \le v\} \in \Psi$ not coinciding
with $(x,y),$ $(x,z)$ or $(y,z),$ and the fiber is the direct product of
an open tetrahedron with vertices called
\begin{equation}
\label{vert}
(x,y), (x,z), (y,z), (u,v)
\end{equation}

\noindent
and an affine subspace of
codimension $3n$ in $K$ (consisting of all maps $\phi: \R^1 \to \R^n$
such that $\phi(x)=\phi(y)=\phi(z)$  and $\phi(u)=\phi(v)$).

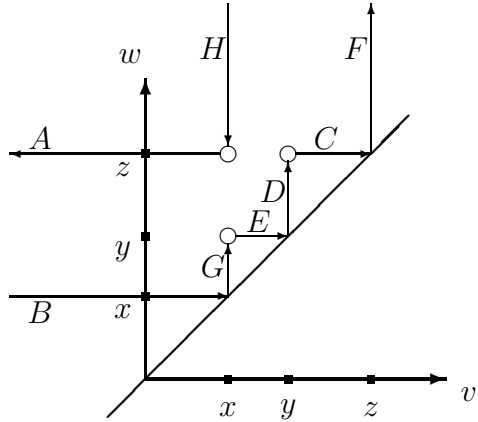
\begin{figure}
\unitlength=1.00mm
\special{em:linewidth 0.4pt}
\linethickness{0.4pt}
\begin{picture}(62.00,60.00)
\thicklines
\put(14.00,1.00){\line(1,1){40.00}}
\put(19.00,6.00){\vector(1,0){40.00}}
\put(19.00,6.00){\vector(0,1){40.00}}
\thinlines
\put(17.00,49.00){\makebox(0,0)[cc]{$w$}}
\put(62.00,4.00){\makebox(0,0)[cc]{$v$}}
\put(29.50,5.50){\rule{1.00\unitlength}{1.00\unitlength}}
\put(37.50,5.50){\rule{1.00\unitlength}{1.00\unitlength}}
\put(38.00,2.00){\makebox(0,0)[cc]{$y$}}
\put(30.00,2.00){\makebox(0,0)[cc]{$x$}}
\put(16.00,23.00){\makebox(0,0)[cc]{$y$}}
\put(18.50,24.50){\rule{1.00\unitlength}{1.00\unitlength}}
\put(18.50,35.50){\rule{1.00\unitlength}{1.00\unitlength}}
\put(16.00,34.00){\makebox(0,0)[cc]{$z$}}
\put(30.00,25.00){\circle{2.00}}
\put(30.00,36.00){\circle{2.00}}
\put(38.00,36.00){\circle{2.00}}
\put(18.50,16.50){\rule{1.00\unitlength}{1.00\unitlength}}
\put(16.00,15.00){\makebox(0,0)[cc]{$x$}}
\put(48.50,5.50){\rule{1.00\unitlength}{1.00\unitlength}}
\put(49.00,2.00){\makebox(0,0)[cc]{$z$}}
\put(1.00,17.00){\vector(1,0){29.00}}
\put(30.00,17.00){\vector(0,1){7.00}}
\put(31.00,25.00){\vector(1,0){7.00}}
\put(38.00,25.00){\vector(0,1){10.00}}
\put(39.00,36.00){\vector(1,0){10.00}}
\put(49.00,36.00){\vector(0,1){20.00}}
\put(30.00,56.00){\vector(0,-1){19.00}}
\put(29.00,36.00){\vector(-1,0){28.00}}
\put(5.00,38.00){\makebox(0,0)[cc]{$A$}}
\put(5.00,15.00){\makebox(0,0)[cc]{$B$}}
\put(28.00,21.00){\makebox(0,0)[cc]{$G$}}
\put(34.00,27.00){\makebox(0,0)[cc]{$E$}}
\put(36.00,31.00){\makebox(0,0)[cc]{$D$}}
\put(43.00,38.00){\makebox(0,0)[cc]{$C$}}
\put(47.00,50.00){\makebox(0,0)[cc]{$F$}}
\put(28.00,50.00){\makebox(0,0)[cc]{$H$}}
\end{picture}
\caption{Boundary of the chain (\protect\ref{first})}
\label{scal1}
\end{figure}

For an arbitrary basepoint $(x,y,z)$ of the former bundle, consider the
intersection set of the fiber over it with the parts of the boundary of
the cycle (\ref{first}), corresponding to 8 graphs
(\ref{aaa3}), (\ref{aaa4}). It is easy to calculate that these
intersections are complete pre-images of the latter
bundle over the open segments in $\Psi \sm ((x,y) \cup (x,z) \cup (y,z))$
labeled in Fig.~\ref{scal1} by characters $A, \ldots, H$, corresponding to
the notation of these graphs in pictures (\ref{aaa3}) and (\ref{aaa4}).

For instance, for any $(4)$-configuration $(r<s<t<u) \in C((4)) = B(\R^1,4)$
the part of the boundary of the fiber in $\Phi_2 \sm \Phi_1$ over this
configuration, corresponding to the graph $E$ in (\ref{aaa4}), belongs
to the fiber of the former bundle in $\Phi_1 \sm \Phi_0$ over the
point $(x<y<z) = (r<t<u)$, and coordinates $(v,w)$ in the base
of the latter bundle satisfy the relations $w=t\equiv y$, $v=s \in (x,y)$.

The part (E) of the boundary of the generator (\ref{first}) of
$\bar H_{\omega-3n+8}(\Phi_2 \sm \Phi_1)$ coincides thus with the complete
pre-image in $\Phi_1 \sm \Phi_0$ of the 4-dimensional submanifold in $C(3,2)$
formed by {\em all} configurations $((x<y<z),(v,w))$ such that $w=y, x<v<y$.

To calculate the coefficients, with which this and other similar pre-images,
corresponding to other letters $A, \ldots, H$,
appear in the {\em algebraic} boundary, we need to fix their orientations.
Again, any such orientation consist of orientations of three objects:

a) the 4-dimensional base manifold in the space of
all configurations $(x,y,z,u,v)$,

b) the tetrahedron (\ref{vert}),
and

c) for any point of this 4-dimensional manifold, the
subspace in $\R^\omega$, consisting of maps $\phi: \R^1 \to \R^n$,
such that $\phi(x)=\phi(y)=\phi(z)$ and $\phi(u)=\phi(v).$

We fix these orientations as follows.

The manifold in the configuration space will be oriented by the pair
of orientations, the first of which is lifted from the
orientation $dx \wedge dy \wedge dz$ of the base of the former bundle, and the
second is given by increase of the unique ``free'' number $u$ or $v$, i.e.
all horizontal segments in Fig.~\ref{scal1} should be oriented to the
right, and all vertical segments should be oriented to the top.

The tetrahedron (\ref{vert}) over any point
$((x<y<z),(v<w))$ of the space $C(3,2)$
will be oriented by the order (\ref{vert}) of its vertices.

The (co)orientation of the fiber $\R^{\omega-3n}$ over the same point
$((x<y<z),(v<w))$ will be specified by the skew form
\begin{eqnarray}
\label{frame2}
\quad d(\xi_1(y)-\xi_1(x)) \wedge d(\xi_2(y)-\xi_2(x)) \wedge \ldots \wedge
d(\xi_n(y)-\xi_n(x)) \wedge \nonumber \\
\wedge \ d(\xi_1(z)-\xi_1(y)) \wedge d(\xi_2(z)-\xi_2(y)) \wedge \ldots \wedge
d(\xi_n(z)-\xi_n(y)) \wedge \nonumber \\
\wedge \ d(\xi_1(w)-\xi_1(v)) \wedge d(\xi_2(w)-\xi_2(v)) \wedge \ldots \wedge
d(\xi_n(w)-\xi_n(v)) \
\end{eqnarray}
on its normal bundle.

Now, for any of 8 components (\ref{aaa3}), (\ref{aaa4}) of the
boundary of the generator (\ref{first}), we need to compare
these orientations with ones induced from the canonical orientations
of similar objects in $\Phi_2 \sm \Phi_1$. We present here these
comparisons only for one component (E).

{\it Compare orientations of configurations.}
In the base $C((4))$ of $\Phi_2 \sm \Phi_1$ the orientation is
given by $dr \wedge ds \wedge dt \wedge du$, and in $\Phi_1 \sm \Phi_0$
the orientation of the corresponding component of the boundary is
given by $dx \wedge dy \wedge dz \wedge dv$. For the
component $E$ we have $x=r, y=t, z=u, v=s$, hence these orientations coincide.

{\it Compare orientations of tetrahedra.} The lexicographic order of
vertices of the graph $E$ in (\ref{aaa4}) is as follows:
$((13), (14), (23), (34))$, or, in notation used in (\ref{vert}),
$((x,y), (x,z), (u,v), (y,z))$. This ordering is of opposite sign with
(\ref{vert}).

{\it Compare coorientations of subspaces $\R^{\omega-3n}$.} If $n$ is even,
then all our coorientations of all such fibers canonically coincide.
For odd $n$ it is not more so, e.g. it is easy to check that for the
component $E$ of the boundary orientations (\ref{frame})
and (\ref{frame2}) are opposite.
\medskip

Doing the same calculations for all other components $A, \ldots, H,$
we obtain Table 1, in which the second line is the sign indicated
in (\ref{aaa3}) or (\ref{aaa4}), and the next three lines are results of
comparisons or orientations of tetrahedra, configuration spaces
and subspaces in $K$ (in the case of odd $n$) respectively.
The 6-th line contains the final sign, with which the corresponding
oriented component, depicted by a segment in Fig.~\ref{scal1},
participates in the boundary of the generator (\ref{first}) in
the case of odd $n$; this sign is just the product of previous
four signs. The last line indicates a similar sign in the case of even $n$;
it is equal to the product of three first signs in the column.

\begin{table}
\begin{center}
\begin{tabular}{|l|c|c|c|c|c|c|c|c|}
\hline
& A & B & C & D & E & F & G & H \\
\hline
Coefficient & -- & + & -- & -- & + & + & -- & + \\
Compare $\Delta$ & -- & -- & + & -- & -- & + & -- & -- \\
Configurations & -- & -- & -- & -- & + & + & + & + \\
Subspaces & + & + & + & -- & -- & + & + & + \\
\hline
Odd $n$ & -- & + & + & + & + & + & + & -- \\
\hline
Even $n$ & -- & + & + & -- & -- & + & + & -- \\
\hline
\end{tabular}
\end{center}
\caption{}
\end{table}

\noindent
\begin{table}
\begin{tabular}{|l|c|c|c|c|c|c|c|c|c|c|c|c|c|c|c|c|c|}
\hline
& A & B & C & D & E & F & G & H & & A & B & C & D & E & F & G & H \\
\hline
Coefficient & -- & + & + & -- & + & -- & -- & + &
& + & -- & + & -- & + & -- & + & -- \\
Compare $\Delta$ & -- & + & -- & + & -- & -- & + & + &
& -- & + & + & -- & -- & -- & -- & -- \\
Configs & + & -- & + & -- & -- & -- & + & + &
& -- & + & -- & + & -- & --& + & + \\
Subspaces & -- & + & + & -- & + & + & + & + &
& -- & + & -- & + & + & + & -- & -- \\
\hline
Odd $n$ & -- & -- & -- & -- & + & -- & -- & + &
& -- & -- & + & + & + & -- & + & -- \\
Even $n$ & + & -- & -- & + & + & -- & -- & + &
& + & -- & -- & + & + & -- & -- & + \\
\hline
\end{tabular}
\caption{}
\end{table}

It follows from these calculations, that in the case of odd $n$ the boundary
in $\Phi_1\sm \Phi_0$ of the element (\ref{first}) is homologous to zero.
Indeed, we can make all signs in the 6-th line equal to + , if defining
the orientation of the 4-submanifold in $C(3,2)$ we orient segments
A and H as shown in Fig.~\ref{scal1}. Therefore our boundary chain
coincides with the boundary of the complete pre-image under the
projection $S(3,2) \to C(3,2)$ of the domain, consisting of all configurations
$((x<y<z),(v<w))$, such that the point $(v,w)$ belongs to the domain in
$\R^2,$ bounded from left and above by segments $H$ and $A$ (see
Fig.~\ref{scal1}), and from right and below by the union of segments
$B,G, E, D, C$ and $F$, depending on $x,y$ and $z$.

In the case of even $n$ the boundary of the element (\ref{first}) is
not a cycle in $\Phi_1 \sm \Phi_0$. However, making all the same calculations
for two other generators
\unitlength=1.00mm
\special{em:linewidth 0.4pt}
\linethickness{0.4pt}
\begin{picture}(28.00,7.00)(2,2)
\put(3.00,6.00){\makebox(0,0)[cc]{1}}
\put(3.00,2.00){\makebox(0,0)[cc]{2}}
\put(13.00,2.00){\makebox(0,0)[cc]{3}}
\put(13.00,6.00){\makebox(0,0)[cc]{4}}
\put(29.00,6.00){\makebox(0,0)[cc]{4}}
\put(29.00,2.00){\makebox(0,0)[cc]{3}}
\put(19.00,2.00){\makebox(0,0)[cc]{2}}
\put(19.00,6.00){\makebox(0,0)[cc]{1}}
\put(16.00,4.00){\makebox(0,0)[cc]{{\large $-$}}}
\put(5.00,1.00){\line(0,1){6.00}}
\put(5.00,7.00){\line(1,-1){6.00}}
\put(11.00,7.00){\line(-1,0){6.00}}
\put(5.00,1.00){\line(1,1){6.00}}
\put(11.00,7.00){\line(0,-1){6.00}}
\put(21.00,7.00){\line(0,-1){6.00}}
\put(21.00,1.00){\line(1,1){6.00}}
\put(27.00,1.00){\line(-1,0){6.00}}
\put(21.00,7.00){\line(1,-1){6.00}}
\put(27.00,1.00){\line(0,1){6.00}}
\end{picture}
and
\unitlength=1.00mm
\special{em:linewidth 0.4pt}
\linethickness{0.4pt}
\begin{picture}(28.00,7.00)(2,2)
\put(3.00,6.00){\makebox(0,0)[cc]{1}}
\put(3.00,2.00){\makebox(0,0)[cc]{2}}
\put(13.00,2.00){\makebox(0,0)[cc]{3}}
\put(13.00,6.00){\makebox(0,0)[cc]{4}}
\put(29.00,6.00){\makebox(0,0)[cc]{4}}
\put(29.00,2.00){\makebox(0,0)[cc]{3}}
\put(19.00,2.00){\makebox(0,0)[cc]{2}}
\put(19.00,6.00){\makebox(0,0)[cc]{1}}
\put(16.00,4.00){\makebox(0,0)[cc]{{\large $-$}}}
\put(11.00,1.00){\line(-1,0){6.00}}
\put(5.00,1.00){\line(1,1){6.00}}
\put(5.00,7.00){\line(0,-1){6.00}}
\put(11.00,1.00){\line(-1,1){6.00}}
\put(5.00,7.00){\line(1,0){6.00}}
\put(27.00,7.00){\line(0,-1){6.00}}
\put(21.00,1.00){\line(1,0){6.00}}
\put(27.00,1.00){\line(-1,1){6.00}}
\put(21.00,7.00){\line(1,0){6.00}}
\put(27.00,7.00){\line(-1,-1){6.00}}
\end{picture}
of the group $\bar H_4(\Delta^2(4))$, we obtain that the {\it difference}
of corresponding chains in $\Phi_2 \sm \Phi_1$ has in
$\Phi_1 \sm \Phi_0$ boundary homologous to zero.
These calculations are represented in Fig.~\ref{aaa63} and  Table 2.
\medskip

Finally, we get the following statement.
\medskip

{\sc Proposition 9.} {\it For any $n$, the group
$\bar H_{\omega-3n+8}(\sigma_3 \sm \sigma_2)$ is isomorphic to $\Z$.
Its basis element coincides in $\Phi_2 \sm\Phi_1$ with the sum
(if $n$ is odd) or difference (for even $n$) of
elements, corresponding to generators
\unitlength=1.00mm
\special{em:linewidth 0.4pt}
\linethickness{0.4pt}
\begin{picture}(28.00,7.00)(2,2)
\put(3.00,6.00){\makebox(0,0)[cc]{1}}
\put(3.00,2.00){\makebox(0,0)[cc]{2}}
\put(13.00,2.00){\makebox(0,0)[cc]{3}}
\put(13.00,6.00){\makebox(0,0)[cc]{4}}
\put(29.00,6.00){\makebox(0,0)[cc]{4}}
\put(29.00,2.00){\makebox(0,0)[cc]{3}}
\put(19.00,2.00){\makebox(0,0)[cc]{2}}
\put(19.00,6.00){\makebox(0,0)[cc]{1}}
\put(16.00,4.00){\makebox(0,0)[cc]{{\large $-$}}}
\put(5.00,1.00){\line(0,1){6.00}}
\put(5.00,7.00){\line(1,-1){6.00}}
\put(11.00,7.00){\line(-1,0){6.00}}
\put(5.00,1.00){\line(1,1){6.00}}
\put(11.00,7.00){\line(0,-1){6.00}}
\put(21.00,7.00){\line(0,-1){6.00}}
\put(21.00,1.00){\line(1,1){6.00}}
\put(27.00,1.00){\line(-1,0){6.00}}
\put(21.00,7.00){\line(1,-1){6.00}}
\put(27.00,1.00){\line(0,1){6.00}}
\end{picture}
and
\unitlength=1.00mm
\special{em:linewidth 0.4pt}
\linethickness{0.4pt}
\begin{picture}(28.00,7.00)(2,2)
\put(3.00,6.00){\makebox(0,0)[cc]{1}}
\put(3.00,2.00){\makebox(0,0)[cc]{2}}
\put(13.00,2.00){\makebox(0,0)[cc]{3}}
\put(13.00,6.00){\makebox(0,0)[cc]{4}}
\put(29.00,6.00){\makebox(0,0)[cc]{4}}
\put(29.00,2.00){\makebox(0,0)[cc]{3}}
\put(19.00,2.00){\makebox(0,0)[cc]{2}}
\put(19.00,6.00){\makebox(0,0)[cc]{1}}
\put(16.00,4.00){\makebox(0,0)[cc]{{\large $-$}}}
\put(11.00,1.00){\line(-1,0){6.00}}
\put(5.00,1.00){\line(1,1){6.00}}
\put(5.00,7.00){\line(0,-1){6.00}}
\put(11.00,1.00){\line(-1,1){6.00}}
\put(5.00,7.00){\line(1,0){6.00}}
\put(27.00,7.00){\line(0,-1){6.00}}
\put(21.00,1.00){\line(1,0){6.00}}
\put(27.00,1.00){\line(-1,1){6.00}}
\put(21.00,7.00){\line(1,0){6.00}}
\put(27.00,7.00){\line(-1,-1){6.00}}
\end{picture}
of the group $\bar H_4(\Delta^2(4))$.

For any $n$, the image of the operator
$d^1: \bar H_{\omega-3n+8}(\Phi_2 \sm \Phi_1) \to
\bar H_{\omega-3n+7}(\Phi_1 \sm \Phi_0)$ is isomorphic to $\Z$
and the quotient group
$\bar H_{\omega-3n+7}(\Phi_1 \sm \Phi_0)/
d^1( \bar H_{\omega-3n+8}(\Phi_2 \sm \Phi_1))$ is isomorphic to $\Z^2$.}
\quad $\Box$

\begin{figure}
\unitlength=1.00mm
\special{em:linewidth 0.4pt}
\linethickness{0.4pt}
\begin{picture}(134.00,65.00)
\thicklines
\put(12.00,10.00){\line(1,1){45.00}}
\put(17.00,15.00){\vector(1,0){40.00}}
\put(17.00,15.00){\vector(0,1){40.00}}
\thinlines
\put(15.00,58.00){\makebox(0,0)[cc]{$w$}}
\put(61.00,12.00){\makebox(0,0)[cc]{$v$}}
\put(27.50,14.50){\rule{1.00\unitlength}{1.00\unitlength}}
\put(35.50,14.50){\rule{1.00\unitlength}{1.00\unitlength}}
\put(36.00,11.00){\makebox(0,0)[cc]{$y$}}
\put(28.00,11.00){\makebox(0,0)[cc]{$x$}}
\put(14.00,32.00){\makebox(0,0)[cc]{$y$}}
\put(16.50,33.50){\rule{1.00\unitlength}{1.00\unitlength}}
\put(16.50,44.50){\rule{1.00\unitlength}{1.00\unitlength}}
\put(14.00,43.00){\makebox(0,0)[cc]{$z$}}
\put(28.00,34.00){\circle{2.00}}
\put(28.00,45.00){\circle{2.00}}
\put(36.00,45.00){\circle{2.00}}
\put(16.50,25.50){\rule{1.00\unitlength}{1.00\unitlength}}
\put(14.00,24.00){\makebox(0,0)[cc]{$x$}}
\put(46.50,14.50){\rule{1.00\unitlength}{1.00\unitlength}}
\put(47.00,11.00){\makebox(0,0)[cc]{$z$}}
\put(31.5,48){\makebox(0,0)[cc]{$A$}}
\put(41,48){\makebox(0,0)[cc]{$B$}}
\put(31,31.5){\makebox(0,0)[cc]{$C$}}
\put(31,39){\makebox(0,0)[cc]{$D$}}
\put(7,36){\makebox(0,0)[cc]{$E$}}
\put(7,29){\makebox(0,0)[cc]{$F$}}
\put(50,55){\makebox(0,0)[cc]{$G$}}
\put(32,55){\makebox(0,0)[cc]{$H$}}
\put(28.00,26.00){\vector(-1,0){24.00}}
\put(28.00,33.00){\vector(0,-1){7.00}}
\put(29.00,45.00){\vector(1,0){6.00}}
\put(28.00,35.00){\vector(0,1){9.00}}
\put(47.00,45.00){\vector(-1,0){10.00}}
\put(47.00,65.00){\vector(0,-1){20.00}}
\put(36.00,46.00){\vector(0,1){19.00}}
\put(4.00,34.00){\vector(1,0){23.00}}
\thicklines
\put(85.00,10.00){\line(1,1){45.00}}
\put(90.00,15.00){\vector(1,0){40.00}}
\put(90.00,15.00){\vector(0,1){40.00}}
\thinlines
\put(88.00,58.00){\makebox(0,0)[cc]{$w$}}
\put(134.00,13.00){\makebox(0,0)[cc]{$v$}}
\put(100.50,14.50){\rule{1.00\unitlength}{1.00\unitlength}}
\put(108.50,14.50){\rule{1.00\unitlength}{1.00\unitlength}}
\put(109.00,11.00){\makebox(0,0)[cc]{$y$}}
\put(101.00,11.00){\makebox(0,0)[cc]{$x$}}
\put(87.00,32.00){\makebox(0,0)[cc]{$y$}}
\put(89.50,33.50){\rule{1.00\unitlength}{1.00\unitlength}}
\put(89.50,44.50){\rule{1.00\unitlength}{1.00\unitlength}}
\put(87.00,43.00){\makebox(0,0)[cc]{$z$}}
\put(101.00,34.00){\circle{2.00}}
\put(101.00,45.00){\circle{2.00}}
\put(109.00,45.00){\circle{2.00}}
\put(89.50,25.50){\rule{1.00\unitlength}{1.00\unitlength}}
\put(87.00,24.00){\makebox(0,0)[cc]{$x$}}
\put(119.50,14.50){\rule{1.00\unitlength}{1.00\unitlength}}
\put(120.00,11.00){\makebox(0,0)[cc]{$z$}}
\put(112,40){\makebox(0,0)[cc]{$A$}}
\put(112,55){\makebox(0,0)[cc]{$B$}}
\put(98,39){\makebox(0,0)[cc]{$C$}}
\put(98,55){\makebox(0,0)[cc]{$D$}}
\put(81,48){\makebox(0,0)[cc]{$E$}}
\put(81,29){\makebox(0,0)[cc]{$F$}}
\put(105,48){\makebox(0,0)[cc]{$G$}}
\put(105,36){\makebox(0,0)[cc]{$H$}}
\put(100.00,34.00){\vector(-1,0){23.00}}
\put(108.00,45.00){\vector(-1,0){6.00}}
\put(101.00,44.00){\vector(0,-1){9.00}}
\put(109.00,65.00){\vector(0,-1){19.00}}
\put(77.00,45.00){\vector(1,0){23.00}}
\put(101.00,46.00){\vector(0,1){19.00}}
\put(106.00,3.00){\makebox(0,0)[cc]{{\large $b$}}}
\put(33.00,3.00){\makebox(0,0)[cc]{{\large $a$}}}
\put(102.00,34.00){\vector(1,0){7.00}}
\put(109.00,34.00){\vector(0,1){10.00}}
\end{picture}
\caption{Boundaries of non-invariant 4-chains}
\label{aaa63}
\end{figure}
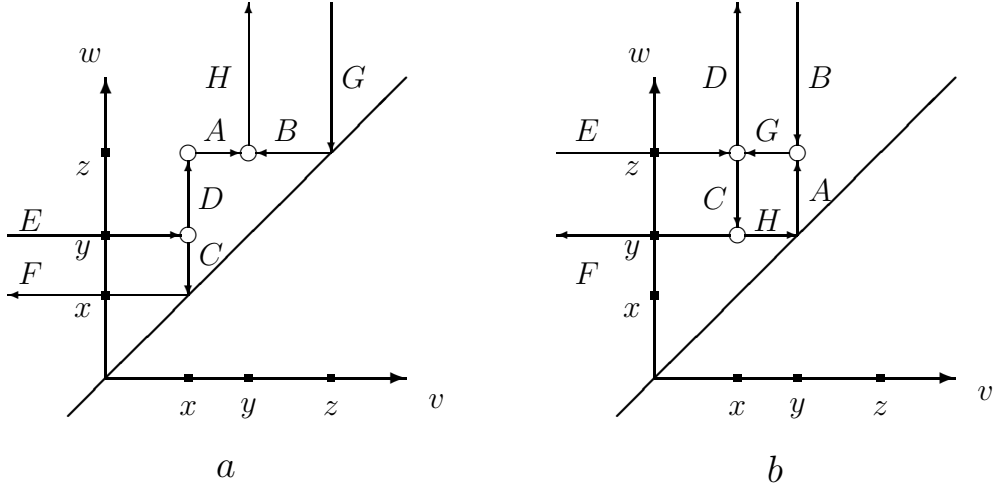

\subsection{The differential
$d^1: \bar H_{\omega-3n+7}(\Phi_1 \sm \Phi_0) \to
\bar H_{\omega-3n+6}(\Phi_0)$.} $\ $

\noindent
{\sc Proposition 10.} {\it The boundary operator
$d^1$ of our spectral sequence maps any of three generators of the group
$\bar H_{\omega-3n+7}(\Phi_1 \sm \Phi_0)$, indicated by
rays in Fig.~\ref{aaa5},
into a generator of the group $\bar H_{\omega-3n+6}(\Phi_0)$.}
\medskip

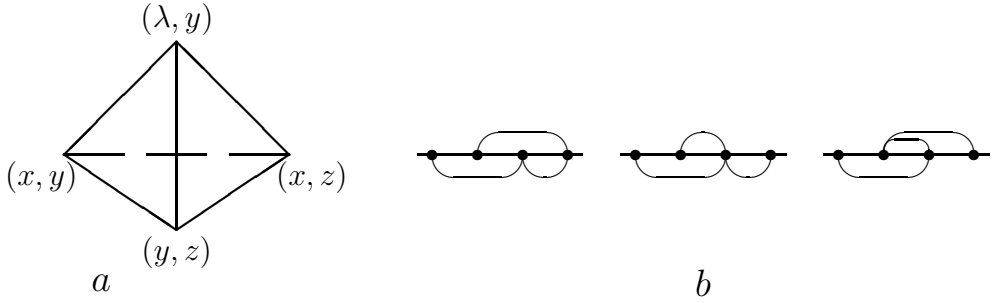
\begin{figure}
\unitlength=1.00mm
\special{em:linewidth 0.4pt}
\linethickness{0.4pt}
\begin{picture}(128.00,38.00)
\put(20.00,7.00){\makebox(0,0)[cc]{$(y,z)$}}
\put(20.00,38.00){\makebox(0,0)[cc]{$(\lambda,y)$}}
\put(38.00,17.00){\makebox(0,0)[cc]{$(x,z)$}}
\put(2.00,17.00){\makebox(0,0)[cc]{$(x,y)$}}
\put(72.00,20.00){\circle*{1.33}}
\put(66.00,20.00){\circle*{1.33}}
\put(60.00,20.00){\circle*{1.33}}
\put(54.00,20.00){\circle*{1.33}}
\put(60.00,20.00){\oval(12.00,6.00)[b]}
\put(69.00,20.00){\oval(6.00,6.00)[b]}
\put(66.00,20.00){\oval(12.00,6.00)[t]}
\put(99.00,20.00){\circle*{1.33}}
\put(93.00,20.00){\circle*{1.33}}
\put(87.00,20.00){\circle*{1.33}}
\put(81.00,20.00){\circle*{1.33}}
\put(126.00,20.00){\circle*{1.33}}
\put(120.00,20.00){\circle*{1.33}}
\put(114.00,20.00){\circle*{1.33}}
\put(108.00,20.00){\circle*{1.33}}
\put(87.00,20.00){\oval(12.00,6.00)[b]}
\put(96.00,20.00){\oval(6.00,6.00)[b]}
\put(90.00,20.00){\oval(6.00,6.00)[t]}
\put(114.00,20.00){\oval(12.00,6.00)[b]}
\put(120.00,20.00){\oval(12.00,6.00)[t]}
\put(117.00,20.00){\oval(6.00,4.00)[t]}
\put(90.00,3.00){\makebox(0,0)[cc]{{\large $b$}}}
\put(10.00,3.00){\makebox(0,0)[cc]{{\large $a$}}}
\thicklines
\put(20.00,10.00){\line(3,2){15.00}}
\put(35.00,20.00){\line(-1,1){15.00}}
\put(20.00,35.00){\line(-1,-1){15.00}}
\put(5.00,20.00){\line(3,-2){15.00}}
\put(20.00,10.00){\line(0,1){25.00}}
\thinlines
\put(5.00,20.00){\line(1,0){8.00}}
\put(16.00,20.00){\line(1,0){8.00}}
\put(27.00,20.00){\line(1,0){8.00}}
\put(52.00,20.00){\line(1,0){22.00}}
\put(79.00,20.00){\line(1,0){22.00}}
\put(106.00,20.00){\line(1,0){22.00}}
\end{picture}
\caption{The boundary of the first generator of
$E^1_{1,\omega-3n+6}$}
\label{tetr}
\end{figure}

We prove this statement for one generator of the former group,
depicted by the horizontal ray in Fig.~\ref{aaa5}, consisting of points
$(\lambda, y)$ with $\lambda<x.$ The chain, representing this generator,
is the fiber bundle, whose base is the configuration space
$B(\R^1,4),$ and the fiber over its point
$\{\lambda<x<y<z\}$ is the direct product of a tetrahedron
(whose vertices are called $(\lambda,y),$
$(x,y),$ $(x,z)$ and $(y,z)$, see Fig.~\ref{tetr}a) and an affine space
of dimension $\omega-3n$ (consisting of all maps $\phi:\R^1 \to \R^n$
such that $\phi(\lambda)=\phi(x)=\phi(y)=\phi(z)$).

One face of the simplex from Fig.~\ref{tetr}a)
(opposite to the vertex $(\lambda,y)$)
is a non-connected graph, and the corresponding stratum
belongs to $\sigma_2$. The strata, swept out by three other faces
(opposite to vertices $(x,y),$ $(x,z)$ and $(y,z)$ respectively),
belong to the term $\Phi_0 \subset \sigma_3 \sm \sigma_2$, namely, they
are complete pre-images in $S(2,2,2)$ of 4-dimensional cycles
in $C(2,2,2),$ consisting of all triples of the form

\begin{equation}
\{(\lambda,y), (x,z), (y,z)\}, \quad
\{(\lambda,y), (x,y), (y,z)\}, \quad
\{(\lambda,y), (x,y), (x,z)\}
\label{faces}
\end{equation}
respectively with all possible
$\lambda<x<y<z$, see Fig.~\ref{tetr}b).
By the Thom isomorphism for the bundle $S(2,2,2) \to C(2,2,2),$
we need only to prove that the union of these three cycles in
$C(2,2,2)$ defines a generator of the group
$\bar H_4(C(2,2,2)) \equiv \bar H_4(B(\Psi,3),{\bf T})$ (if $n$ is odd) or
$\bar H_4(C(2,2,2),\pm \Z) \equiv \bar H_4(B(\Psi,3),{\bf T};\pm \Z)$
(if $n$ is even);
or equivalently, that its boundary in $\bar H_3({\bf T})$
is the fundamental cycle of ${\bf T} \sim B(\R^1,3).$
But among three manifolds (\ref{faces}) only the first approaches
${\bf T}$ and its boundary obviously coincides with ${\bf T}.$

The calculation for two other generators of the group
$\bar H_{\omega-3n+7}(\Phi_1 \sm \Phi_0)$ is exactly the same.
\quad $\Box$ \medskip

Theorem 5 is a direct corollary of Lemma 6 and Propositions
9 and 10. \quad $\Box$

\section{Cohomology classes of compact knots in $\R^n$.}

\subsection{Classes of order 1}

It is well-known that there are no first-order knot invariants in $\R^3$,
see \cite{V1}. However, the subgroup
$F_{1,\Z_2}^* \subset H^*(\K \sm \Sigma, \Z_2)$
of all $\Z_2$-valued first-order cohomology classes of the space of
compact knots is non-trivial. For instance, the generator of the
group $F^1_{1,\Z_2} \subset H^1(\K \sm \Sigma, \Z_2)$ proves that
already the component of unknots in $\R^3$ is not simple-connected.

More generally, the following statements hold.
\medskip

{\sc Theorem 6.} {\it A. For any $n \ge 3$, the subgroup
$F_{1,\Z_2}^* \subset H^*(\K \sm \Sigma, \Z_2)$ of
first-order cohomology classes of the space of knots in $\R^n$
contains exactly two non-trivial homogeneous
components $F_{1,\Z_2}^{n-2} \sim F_{1,\Z_2}^{n-1} \sim \Z_2$.

B. If $n$ is even, then both these cohomology classes give rise to
integer cohomology classes, i.e. $F_{1,\Z}^{n-2} \sim F_{1,\Z}^{n-1} \sim \Z$,
and there are no other non-trivial integer cohomology
groups $F_{1,\Z}^d$, $d \ne n-2, n-1$.

C. If $n$ is odd, then the generator of the group
$F_{1,\Z_2}^{n-1}$ is equal to the first Steenrod operation
of the generator of
$F_{1,\Z_2}^{n-2}$.

D. The generator of the group
$F_{1,\Z_2}^{n-2}$  can be defined as
the linking number with the $\Z_2$-fundamental cycle
of the variety $\Gamma \subset \Sigma$, formed by all
maps $\phi: S^1 \to \R^n$, gluing together some two {\it opposite} points of
$S^1$; the generator of the group $F_{1,\Z_2}^{n-1}$ is the linking number with
the $\Z_2$-fundamental cycle of the subvariety $\Gamma! \subset \Gamma$,
formed by all maps $\phi: S^1 \to \R^n$, gluing together some two
{\bf fixed} opposite points, say, the points $0$ and $\pi.$
Moreover, if $n$ is even, then these two varieties are orientable, and
the groups
$F_{1,\Z}^{n-2}$, $F_{1,\Z}^{n-1}$ are generated by the linking numbers
with the corresponding $\Z$-fundamental cycles.

E. If $n=3,$ then the cycles, generating the groups
$F_{1,\Z_2}^{1}$ and $F_{1,\Z_2}^{2}$ are non-trivial
already in the restriction on the component of the unknot in
$\R^3.$} \medskip

{\it Proof.} The term $\sigma_1$ of the main filtration of
$\sigma$ is the space of an
affine fiber bundle of dimension $\omega-n$, whose base
is the configuration space $\Psi$; in the case of compact
knots this base space is diffeomorphic to the closed
M\"obius band. It is easy to check that this affine bundle is
orientable if and only if $n$ is even. Thus the
term $E_1^{-1,q}(\Z_2)$ of the main spectral sequence with coefficients
in $\Z_2$ is isomorphic to $\Z_2$ for $q= n-1$ or $n$ and is trivial
for all other $q$. Moreover, if $n$ is even, then
also the terms $E_1^{-1,q}(\Z)$  are isomorphic to
$\Z$ if $q= n-1$ or $n$ and are trivial for all other $q.$
The basic cycles in $\sigma_1$, generating these groups, are
the manifolds $\tilde \Gamma$ (respectively, $\tilde \Gamma !$),
consisting of all pairs $(\phi,(x,y))$ such that $\phi(x)=\phi(y)$
and $y=x+\pi$ (respectively, $x=0, y=\pi$). Thus the direct images of these
cycles in $\Sigma$ are exactly the fundamental classes of varieties
$\Gamma$, $\Gamma!$, mentioned
in statement D of Theorem 6. We need to prove that the homology
classes in $\Sigma$ of these cycles are
non-trivial, i.e. the linking numbers with them
(which we denote by $\{\Gamma\}$ and
$\{\Gamma!\}$ respectively) are non-trivial cohomology
classes in $\K \sm \Sigma.$

For $n>4$ this follows from the dimensional reasons. Indeed, for any
$n >3$ our spectral sequence converges to the entire group
$H^*(\K \sm \Sigma),$ and by (\ref{ineq}), (\ref{ineq2}) all groups
$E_{\ge 1}^{p,q}$, $p<-1,$ with $p+q< n-1$ are trivial, hence all differentials
$d^r$, acting into the cells $E_r^{-1,q}$, also are trivial. For $n=4$ there is
unique non-trivial group $E_1^{-2,4},$ from which, in principle, there could
act a non-trivial homomorphism $d^1 : E_1^{-2,4} \to E_1^{-1,4} \sim \Z.$
We will see in the next subsection, that the group $E_1^{-2,4}(G)$
is isomorphic to the basic coefficient group $G.$ Thus for
$G=\Z$ and $n=4$ the differential $d^1(E_1^{-2,4})$ is trivial
by the Kontsevich's realization theorem (\ref{ineq3}), and for $G=\Z_2$
the same follows from the functoriality of spectral sequences
under the coefficient homomorphisms. Thus statements A, B, D
of Theorem 6 are proved for all $n>3$.

The statement C for odd $n>3$ will be proved as follows. First, we
prove the equality $Sq^1(\{\Gamma\}) = \{\Gamma!\}$ in
the group $H^*(\K \setminus \Gamma),$ then the desired similar equality
in the cohomology group of the subspace
$\K \setminus \Sigma \subset \K \setminus \Gamma$ will follow
from the naturality of Steenrod operations.

For any $n \ge 3$, the group $H^*(\K \sm \Gamma)$ can be calculated by a
cohomological spectral sequence similar to (but much more easy than)
our main spectral sequence. Namely, it is ``Alexander dual''
(in the sense of formula (\ref{dual}))
to the homological sequence
associated to the standard filtration of the
simplicial resolution $\gamma$ of $\Gamma$. The term
$\gamma_i \sm \gamma_{i-1}$
of this filtration is the space of a fiber bundle with base
$B(S^1,i)$ (where $S^1$ is the space of all diametral
pairs of points of the original circle) and the fiber equal to
the direct product of the open simplex
$\dot \Delta^{i-1} $ and the space $\R^{\omega-ni}$. In particular
$E_1^{p,q}=0$ for $p>0$ or $q+p(n-1)<0.$ For any $n>3$ the group
$H^{n-1}(\K \sm \Gamma,\Z_2)$, containing the element
$Sq^1(\{\Gamma\})$, is one-dimensional and is generated by
the linking number with the cycle $\Gamma!$. Consider any loop
$l \subset \Gamma,$ lying in the set of regular points of $\Gamma$
and such that going along $l$ we permute two points $x, y=x+\pi \in S^1$,
glued together by the corresponding maps.
(Such $l$ exist for any $n \ge 3$, because the
codimension in the manifold $\tilde \Gamma \simeq \gamma_1$ of the preimage of
the set of singular points of $\Gamma$ is equal to $n-1$, and we can realize $l$
as the projection into $\Gamma$ of almost any loop in $\tilde \Gamma$,
permuting the points $x$ and $y$.)
Let $L$ be the ``tube'' around $l$ in $\K \sm \Gamma,$ i.e. the
union of boundaries of small $(n-1)$-dimensional discs
transversal to $\Gamma$ with centers  at the points
of $l$.
The fibration $L \to l$ with
fiber $S^{n-2}$ is non-orientable if $n$ is odd, thus already in the restriction
to $L$ $Sq^1(\{\Gamma\})$ is non-trivial and is equal to the $\Z_2$-fundamental
cocycle of $L$.  The union of these transversal discs spans $L$ in $\K$, and its
$\Z_2$-intersection number with $\Gamma!$ is equal to 1, in particular the
fundamental cycle of $L$ is nontrivial in $H_{n-1}(\K \sm \Gamma,\Z_2)$ for any
$n \ge 3$ and generates this group for $n > 3$. Also, we get that for any odd $n
\ge 3$ in restriction to $L$ $Sq^1(\{\Gamma\})=\{\Gamma!\}$.  Since for $n>3$
the inclusion homomorphism $H^{n-1}(\K \sm \Gamma,\Z_2) \to H^{n-1}(L,\Z_2)$ is
an isomorphism, statement C of Theorem 6 is proved for all odd $n>3$.

\begin{figure}
\unitlength=1.00mm
\special{em:linewidth 0.4pt}
\linethickness{0.4pt}
\begin{picture}(130,48)
\put(12,23){\oval(24,6)}
\put(18,43){\oval(6,6)[t]}
\put(18,43){\oval(6,6)[bl]}
\put(18,41){\oval(6,2)[br]}
\put(30,43){\oval(18,6)[b]}
\put(30,43){\oval(18,6)[tr]}
\put(30,45){\oval(18,2)[tl]}
\put(52,43){\oval(12,6)[t]}
\put(52,43){\oval(12,6)[bl]}
\put(52,41){\oval(12,2)[br]}
\put(64,43){\oval(12,6)[b]}
\put(64,43){\oval(12,6)[tr]}
\put(64,45){\oval(12,2)[tl]}
\put(85,43){\oval(18,6)[t]}
\put(85,43){\oval(18,6)[bl]}
\put(85,41){\oval(18,2)[br]}
\put(97,43){\oval(6,6)[b]}
\put(97,43){\oval(6,6)[tr]}
\put(97,45){\oval(6,2)[tl]}
\put(107,23){\oval(24,6)}
\put(85,3){\oval(18,6)[b]}
\put(85,3){\oval(18,6)[tl]}
\put(85,5){\oval(18,2)[tr]}
\put(97,3){\oval(6,6)[t]}
\put(97,3){\oval(6,6)[br]}
\put(97,1){\oval(6,2)[bl]}
\put(52,3){\oval(12,6)[b]}
\put(52,3){\oval(12,6)[tl]}
\put(52,5){\oval(12,2)[tr]}
\put(64,3){\oval(12,6)[t]}
\put(64,3){\oval(12,6)[br]}
\put(64,1){\oval(12,2)[bl]}
\put(18,3){\oval(6,6)[b]}
\put(18,3){\oval(6,6)[tl]}
\put(18,5){\oval(6,2)[tr]}
\put(30,3){\oval(18,6)[t]}
\put(30,3){\oval(18,6)[br]}
\put(30,1){\oval(18,2)[bl]}
\put(12,20){\circle*{1.33}}
\put(27,40){\circle*{1.33}}
\put(58,43){\circle*{1.33}}
\put(88,46){\circle*{1.33}}
\put(107,26){\circle*{1.33}}
\put(88,6){\circle*{1.33}}
\put(27,0){\circle*{1.33}}
\end{picture}
\caption{Non-trivial 1-cycle in the space of unknots}
\label{ha}
\end{figure}
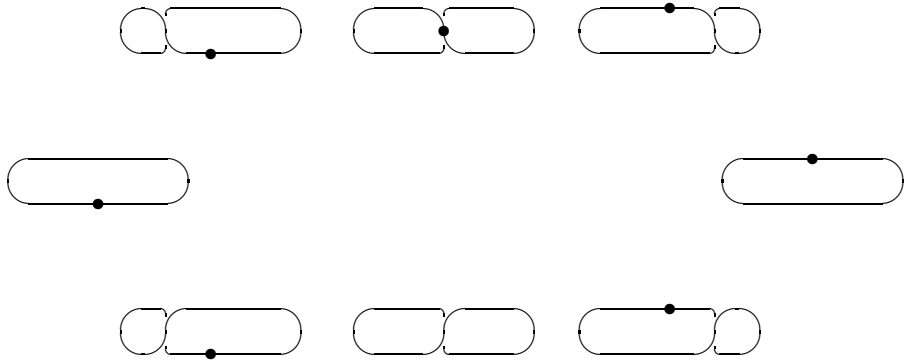

Finally, let be $n=3.$ Consider the loop $\Lambda:S^1 \to (\K \sm \Sigma),$
some whose eight points are shown in Fig.~\ref{ha}. Note that
any two (un)knots of this family, placed in this picture one over the
other, have the same projection to $\R^2$. Let us connect any such two
unknots by a segment in $\K$, along which the projection to $\R^2$ also is
preserved. The union of these segments is a disc in $\K$, spanning
the loop $\Lambda$; it is obvious that the
$\Z_2$-intersection number of this disc with
the variety $\Gamma$ is equal to 1, in particular the class
$\{\Gamma\} \in H^1(\K \sm \Sigma, \Z_2)$ is non-trivial.

It is easy to see that this loop $\Lambda$ is homotopic to the loop
$\Lambda'$, consisting of knots, obtained from the standard
embedding $\phi: S^1 \to \R^2 \subset \R^3$ by rotations by all
angles $\al \in [0,2\pi]$ around any diagonal of $\phi(S^1)$,
and also to the loop $\Lambda''$, consisting of all knots
$\phi_\tau,$ $\tau \in [0,2\pi],$ having the same image as $\phi$
and given by the formula $\phi_\tau(\alpha) = \phi(\alpha+\tau).$

Let us fix some sphere $S^2 \subset \R^3$ and consider the space $GC$
of all naturally pa\-ra\-met\-ri\-zed great circles in it.
This space is obviously diffeomorphic to $SO(3) \sim \R P^3,$
and its group $H_1(GC,\Z_2)$ is generated by our loop $\Lambda'
\sim \Lambda.$ To complete the proof of statements A, C, D for $n=3$
and of statement E, it remains to show that the linking number of
the generator of $H_2(GC,\Z_2)$ with the variety
$\Gamma!$ is equal to 1.

If $n=3,$ then the group $H^2(\K \sm \Gamma),$ containing the element
$Sq^1(\{\Gamma\}),$ is two-dimensional: besides the class
$\{\Gamma!\}$ it contains a basic element of second order, coming
from the cell $E^{-2,4}$ of the canonical spectral sequence calculating
$H^*(\K \sm \Gamma, \Z_2)$. The reduction mod $\gamma_1$ of its dual class
in $\bar H_{\omega-3}(\Gamma)$ is the fundamental cycle of
$\gamma_2 \sm \gamma_1$. As we will see in the next subsection, this
fundamental cycle is homologous to zero in the space
$\sigma_2 \sm \sigma_1 \supset \gamma_2 \sm \gamma_1.$ Thus
the element $Sq^1(\{\Gamma\}) \in H^2(\K \sm \Sigma, \Z_2)$
belongs to the group of elements of order 1, which is generated
by  the class $\{\Gamma!\}.$ On the other hand,
$Sq^1(\{\Gamma\}) \equiv \{\Gamma\}^2$ is non-trivial already in
restriction to the submanifold $GC \subset \K \sm \Sigma$
(because $\{\Gamma\}$ is), thus it coincides with $\{\Gamma!\}.$
Theorem 6 is completely proved. \quad $\Box$

\subsection{Classes of order 2}

{\sc Theorem 7.} {\it For any $n \ge 3,$ the group
$F^*_{2,\Z} / F^*_{1,\Z}$ of order 2 cohomology classes of the space
of compact knots, reduced modulo classes of order 1, contains at most two
non-trivial homogeneous components in dimensions
$2n-6$ and $ 2n-3$. The component in dimension $2n-6$ is always
isomorphic to $\Z$, the component in dimension $2n-3$ is isomorphic
to $\Z$ if $n>3$ and is a cyclic group if $n=3.$
The generator of the nontrivial group in
dimension $2n-3$ has bi-order $(2,0)$, and the generator in dimension
$2n-6$ has bi-order $(2,1)$.} \medskip

{\it Proof.} The auxiliary filtration of the space $\sigma_2 \sm \sigma_1$
consists of two terms $\Phi_0 \subset \Phi_1.$ $\Phi_0$ is the space
of a fiber bundle over $B(\Psi,2)$, whose fiber is the direct product of
a vector subspace of codimension $2n$ in $\K$  and an open interval. It is
easy to calculate that the Borel--Moore homology group
$\bar H_*(\Phi_0)$ of the space of this bundle coincides with that of its
restriction on the subset $B(S^1,2) \subset B(\Psi,2),$ where $S^1$ is the
{\it equator} of the M\"obius band $\Psi$, consisting of all pairs of the
form $(\alpha, \alpha+\pi)$. This subset $B(S^1,2)$ is
homeomorphic to the open M\"obius band. The generator of its
fundamental group $\pi_1(B(S^1,2)) \simeq \Z$ preserves the orientation
of the bundle of $(\omega-2n)$-dimensional subspaces in $\K,$
and destroys the orientations of both the bundle of open intervals and
the base $B(S^1,2)$ itself. Therefore $\bar H_*(\Phi_0)$ coincides
with the homology group of a circle up to the shift of dimensions
by $\omega-2n+2$, namely, it has only two non-trivial homology
groups $\bar H_{\omega-2n+3}$ and $\bar H_{\omega-2n+2},$
which are isomorphic to $\Z$. The first (respectively, second)
of them is generated by the fundamental cycle of the preimage of the set
$B(S^1,2) \subset B(\Psi,2)$  (respectively, by that of the
circle in $B(S^1,2)$ formed by all pairs of pairs of the form
$((\al,\al+\pi),(\al+\pi/2,\al+3\pi/2))$, where $\al$ is defined up to
addition of a multiple of $\pi/2$).

The term $\Phi_1 \sm \Phi_0$ is the space of a fiber bundle over
$B(S^1,3)$, whose fiber is the direct product of a vector
subspace of codimension $2n$  in $\K$ and an open triangle.
Both these bundles (of vector subspaces and triangles) over $B(S^1,3)$ are
orientable, thus $\bar H_*(\Phi_1 \sm \Phi_0) \simeq
\bar H_{*-(\omega-2n+2)}(B(S^1,3)),$ i.e., this group is equal to
$\Z$ in dimensions $\omega-2n+5$ and $\omega-2n+4$
and is trivial in all other dimensions. The
$(\omega-2n+5)$-dimensional component is generated by the
fundamental cycle of $\Phi_1 \sm \Phi_0,$ and the
$(\omega-2n+4)$-dimensional one by the preimage of the
cycle in $B(S^1,3)$, formed by all triples
$(\alpha,\beta,\gamma)\subset S^1$
such that $\alpha+\beta+\gamma \equiv 0 \, ($mod $2\pi).$

The unique nontrivial differential $d^1: \E^1_{1,*} \to \E^1_{0,*}$
of the reversed auxiliary spectral sequence sends the
$(\omega-2n+4)$-dimensional generator of the group
$\bar H_*(\Phi_1 \sm \Phi_0)$ to a multiple of the
$(\omega-2n+3)$-dimensional generator of $\bar H_*(\Phi_0).$
Now we prove that the coefficient in this multiple is equal to $\pm 1.$

By construction, the image of this differential is the fundamental
cycle of the preimage in $\Phi_0 $ of a certain 2-dimensional
cycle $\Delta \subset B(\Psi,2)$. Namely,
this cycle is the space of a fiber bundle
over $S^1,$ whose fiber over the point $\alpha \in S^1$ is homeomorphic to an
open interval and consists of all unordered pairs $\{(\alpha,\beta);
(\alpha, \gamma)\}$  such that $\beta \ne \gamma$ and
$\alpha+\beta+\gamma \equiv 0 \, ($mod $2\pi)$
(although either $\beta$ or $\gamma$ can coincide with $\alpha$).
\medskip

{\sc Lemma 7.} {\it The submanifolds $\Delta$ and $B(S^1,2)$ of $B(\Psi,2)$
are homeomorphic, and there is a proper homotopy between this homeomorphism
$\Delta \to B(S^1,2)$ and the identical embedding of $\Delta$ into
$B(\Psi,2)$.} \medskip

\begin{figure}
\unitlength=1.00mm
\special{em:linewidth 0.4pt}
\linethickness{0.4pt}
\begin{picture}(66.00,23.00)
\put(20.00,13.00){\circle{14.00}}
\put(55.00,13.00){\circle{14.00}}
\put(26.00,10.00){\line(-1,6){1.33}}
\put(56.00,6.00){\line(-1,6){2.33}}
\put(61.00,16.00){\line(-2,-1){12.00}}
\put(26.00,10.00){\line(-2,-1){8.00}}
\put(32.00,13.00){\vector(1,0){11.00}}
\put(18.00,3.00){\makebox(0,0)[cc]{$\gamma$}}
\put(26.00,21.00){\makebox(0,0)[cc]{$\beta$}}
\put(28.00,8.00){\makebox(0,0)[cc]{$\alpha$}}
\put(45.00,8.00){\makebox(0,0)[cc]{$\gamma(1)$}}
\put(56.00,3.00){\makebox(0,0)[cc]{$\alpha'(1)$}}
\put(67.00,16.00){\makebox(0,0)[cc]{$\alpha''(1)$}}
\put(54.00,23.00){\makebox(0,0)[cc]{$\beta(1)$}}
\end{picture}
\caption{The homotopy $\Delta \to B(S^1,2)$}
\label{X}
\end{figure}
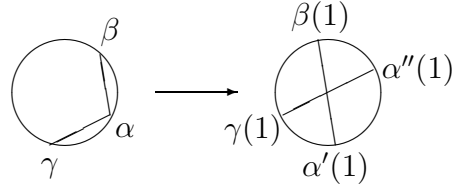

{\it Proof.} We will consider all points $\alpha, \beta$ etc.
as points of the unit circle in the complex line $\C^1.$
Any point $\{(\alpha,\beta);(\alpha,\gamma)\} \in \Delta$
can be deformed to a point of $B(S^1,2)$ by a homotopy $h_t, $ $t \in [0,1]$
(see Fig. \ref{X}), along which the chords $[\alpha'(t),\beta(t)]$
and $[\alpha''(t),\gamma(t)]$ preserve their directions.
The easiest possible ``direct and uniform'' realization of such a
deformation depends continuously on the initial point of $\Delta$;
it is easy to calculate that any point of $B(S^1,2)$  can be obtained by such a
homotopy from exactly one point of $\Delta.$ \quad $\Box$
\medskip

So, the differential $d^1:\E^1_{1,*} \to \E^1_{0,*}$ kills both groups
$\bar H_{\omega-2n+4}(\Phi_1 \sm \Phi_0)$
and $\bar H_{\omega-2n+3}(\Phi_0)$, and
the column $E_1^{-2,*}$ of the main cohomological spectral sequence
contains exactly two non-trivial terms
$E_1^{-2,2n-4} \sim E_1^{-2,2n-1} \sim \Z.$ By the Kontsevich's
realization theorem (see \S \ 4.2), they survive after the action of
all  differentials of the spectral sequence and define certain
cohomology classes of $\K \sm \Sigma$ in dimensions $2n-6$
and $2n-3$ respectively. However, in the case $n=3$ we cannot
be sure that the $(2n-3)$-dimensional class or some its multiple
will not be killed by some unstable cycle in $\Sigma$, not
counted by our spectral sequence. The $(2n-6)$-dimensional
class in this case is a well-known knot invariant, and Theorem 7
is proved.

Note however, that for $n=4$ its statement, concerning the group in
dimension $2n-3,$ depends on a non-published theorem of Kontsevich.
(For $n>5$ the fact that nothing acts into the corresponding cell
$E^{-2,2n-1}$ follows from the dimensional reasons,
see (\ref{ineq2}), and for $n=5$ from the simplest version
(\ref{ineq3}) of the Kontsevich's theorem, whose proof essentially coincides
with that for the classical case $n=3$, see \cite{K}, \cite{BN2}.  \quad $\Box$

\section{Three problems}

{\sc Problem 1.} To establish a direct correspondence between
the algebra of two-connected graphs and that of Chinese Character Diagrams,
see \cite{BN2}. Perhaps this correspondence will give us the most natural
proof of Theorem 2. The starting point of our construction was the
space of all smooth maps $\R^1 \to \R^n.$ Is there some analog
of it behind the Chinese Character Diagrams?
\medskip

{\sc Problem 2.} To present a precise (and economical) description of
generators of groups $\bar H_*(\Delta^2(k))$ The
absence of such a description is now almost unique objective
to the calculation of order 4 cohomology classes of spaces of
non-compact knots.
\medskip

{\sc Problem 3.} The multiplication conjecture of \S \ 5.5.2.

\end{document}